\documentclass{article}
\usepackage{amsmath}
\usepackage{amssymb}
\usepackage{amsthm}
\usepackage{color}
\usepackage{tikz}
\usepackage[colorlinks,citecolor=red,pagebackref=false,hypertexnames=false]{hyperref}

\newtheorem{thm}{Theorem}[section]
\newtheorem{lemma}{Lemma}[section]
\newtheorem{remark}{Remark}[section]

\newtheorem{proposition}{Proposition}[section]
\title{Recovering point sources for the inhomogeneous Helmholtz equation}
\date{}

\author{Gang Bao\footnote{Department of
Mathematics, Zhejiang University, Hangzhou,
Zhejiang, 310027, China. Email: baog@zju.edu.cn} \; Yuantong Liu\footnote{Department of
Mathematics, Zhejiang University, Hangzhou,
Zhejiang, 310027, China. Email: ytliu@zju.edu.cn} \; \textrm{and}
Faouzi Triki\footnote{Laboratoire Jean Kuntzmann, UMR CNRS 5224, 
Universit\'{e} Grenoble-Alpes, 700 Avenue Centrale,
38401 Saint-Martin-d'H\`eres, France. Email: faouzi.triki@univ-grenoble-alpes.fr. 
The work of GB is supported in part by a NSFC Innovative Group Fund (No.11621101).
The work of FT is supported in part by the grant ANR-17-CE40-0029 of the French National Research Agency ANR (project MultiOnde).}}

\begin{document}
\maketitle

\begin{abstract}

    The paper is concerned with an inverse point source problem for the Helmholtz equation. It consists of recovering the locations and amplitudes of a finite number of radiative point sources inside a given inhomogeneous medium from the knowledge of a single boundary measurement. The main result of the paper is a new H\"{o}lder type stability estimate for the inversion under the assumption that the point sources are well separated. The proof of the stability is based on a combination of Carleman estimates and a technique for proving uniqueness of the Cauchy problem.
\end{abstract}

\section{Introduction}
We investigate the issues of uniqueness and stability for recovering multiple point sources for a Helmholtz equation in an inhomogeneous medium from boundary measurements. This inverse problems arise in a diverse set of application areas, including medical imaging \cite{fokas2004unique, hamalainen1993magnetoencephalography, he1998identification}, antenna synthesis \cite{angel1991antenna}, pollution detection \cite{elbadia2005inverse}. Consider the radiated field $u$ at a frequency $k>0$, by a compactly supported source function $S$, that satisfies the following Helmholtz equation

\begin{align}\label{equation: forward problem}
    \begin{cases}
        \Delta u + k^{2}(1+q(x)) u = S \quad \text{in} \quad \mathbb{R}^{3}, \\
        \lim\limits_{r\to +\infty}r(\frac{\partial u}{\partial r}-iku)=0,\quad r=\left\|x \right\|,
    \end{cases}
\end{align}

where $q$ is the medium function. Given the source $S$, the forward problem is to determine the field $u$ while the inverse source problem is to recover $S$ from near-field or far-field measurements of $u$.

For a general source $S$, the inverse source problem is ill-posed. The uniqueness is not guaranteed due to the existence of non-radiating sources \cite{alves2009full, bleistein1977nonuniqueness, bao2010multifrequency}. To overcome this difficulty, further information on the source or different kinds of measurements is needed. Indeed the identification of a general source may be achieved for example by considering multifrequency boundary measurements \cite{bao2010multifrequency, bao2011numerical, acosta2012multifrequency}. Unlike general source functions, point sources are singular and have a lower dimensionality. These properties are sufficient to uniquely determine the point source using a single boundary measurement \cite{ammari2002inverse}. Many stability estimates of H\"{o}lder type have been derived for this inverse problem under the assumption that the background medium is homogeneous \cite{elbadia2012holder, elbadia2011inverse}. The proofs are mainly based on solving an algebraic system obtained by using test functions in the form of products of polynomial functions and plane waves. In \cite{ren2019imaging}, Ren and Zhang established a H\"{o}lder stability estimate on the reconstruction of point sources with respect to smooth changes of a known medium. Their results indicate that if the medium is known up to a small smooth perturbation, the recovered source is close to the real one in a topology induced by a nonconventional given metric. We refer to \cite{bao2011numerical, bao2015recursive, liu2020reconstruction, wang2019fourier} for numerical treatments of the inverse problem. Note that the obtained stability estimates are useful for proving error estimates of the reconstruction algorithms.

This work is a generalization of a stability estimate in Sobolev spaces with negative exponents obtained recently in \cite{bao2021recovering} for a single point source. We are interested in recovering the locations $z_{i}$ and intensities $a_{i}$ of the source $S= \sum_{i=1}^{N} a_{i} \delta_{z_{i}}$ from the knowledge of the field $u$ on the boundary of a given domain $\Omega$ containing the point sources as well as the support of the medium function $q$. Note that since $u$ satisfies the Helmholtz equation with zero right hand side outside $\Omega$, and Sommerfeld radiation conditions, the knowledge of the field on the boundary of $\Omega$ is equivalent to having the Cauchy data of the field.

Let $\Omega$ be a bounded domain in $\mathbb{R}^{3}$ with a $C^{2}$ boundary $\partial \Omega$. Denote $\nu(x)$ the outward unit normal vector at $x\in \partial \Omega$. Assume that $ q\in C^{2}(\mathbb{R}^{3})$ has a compact support in $\Omega$, and let $R_{0}>0$ be large enough such that $\Omega \subset B_{\frac{R_{0}}{2}}$, here we use the notation $B_{R}=\{ x\in \mathbb{R}^{3}: \,\lvert x \rvert \leq R \}$ and $B_{R}(z)=\{ x\in \mathbb{R}^{3}: \,\lvert x-z \rvert \leq R \}$ for $z\in \mathbb{R}^{3},R\in \mathbb{R}$. For given integer $N_{0}>0$, a constant $ \overline{a}>0$, and a constant $\eta>0$ small enough, define the set

\begin{align}
    \mathcal P = \left\{ \left(a_{j}, z_{j}\right)_{j=1,\cdots N} \in \mathbb {R}^{4N}: N\in \mathbb N^{*},\, N\leq N_{0},\, 0 \leq a_{j}<\overline{a}, \, 8\eta<\lvert z_{i}-z_{j} \rvert,\forall i\neq j,\, 8\eta < \textrm{dist}(z_{i}, \partial\Omega)\right\}.
    \label{S}
\end{align}

We further denote $\mathfrak p= (\overline a, \eta, N_{0})$, $\mathcal S$ the set of point sources $S= \sum_{i=1}^{N} a_{i} \delta_{z_{i}}$, with features $\left(a_{j}, z_{j}\right)_{j=1,\cdots N} \in \mathcal P $. This geometrical assumption plays a crucial role in deriving our stability estimate for point sources within the set $\mathcal S$.

Now, we are ready to present our main stability estimate.

Let $S_{1} = \sum_{i=1}^{N_{1}} a_{1,i} \delta_{z_{1,i}}$ and $S_{2}= \sum_{i=1}^{N_{2}} a_{2,i} \delta_{z_{2,i}}$ be two sources in $\mathcal S$. We can show that $\lvert z_{i,j}-z_{i,k} \rvert \geq 8\eta$ $i=1, 2, 1\leq j\neq k \leq N_{i}$ implies that $B_{3\eta}(z_{i,j})$ contains a maximum of two point sources.

We further denote $\mathfrak N \subset \{1, \cdots, N_{1}\}$ the set for which $B_{3\eta}(z_{1,j}), j \in \mathfrak N,$ holds exactly two point sources, and let $\pi: \mathfrak N \rightarrow \{1, \cdots, N_{2}\}$ be a one-to-one mapping such that $z_{2,\pi(j)} \in B_{3\eta}(z_{1,j})$ for all $j\in \mathfrak N$. Let $\mathfrak N_{1} = \{1, \cdots, N_{1}\} \setminus \mathfrak N $, $\mathfrak N_{2} = \{1, \cdots, N_{2}\} \setminus \pi(\mathfrak N). $

\begin{thm}\label{thm: stability for multiple a z} Let $S_{1} = \sum_{i=1}^{N_{1}} a_{1,i} \delta_{z_{1,i}}$ and $S_{2}= \sum_{i=1}^{N_{2}} a_{2,i} \delta_{z_{2,i}}$ be two sources in $\mathcal S$. Let $u_{1},u_{2}$ be the solutions to the inhomogeneous Helmholtz equation \eqref{equation: forward problem} with respectively point sources $S_{1} $ and $S_{2}$, and denote
    $$
        \varepsilon= \lVert u_{1}-u_{2} \rVert_{H^{1}(\partial \Omega)} + \lVert \frac{\partial u_{1}}{\partial\nu} - \frac{\partial u_{2}}{\partial \nu} \rVert_{L^{2}(\partial \Omega)}.
    $$
    If $\varepsilon$ is small enough such that $\varepsilon\leq\varepsilon_{0}$, then there exist positive constants $C_{\mathfrak{a}}, C_{\mathfrak{z}},$ 
    such that
    \begin{align*}
        \lvert a_{1,i}-a_{2,\pi(i)} \rvert
        \leq
        C_{\mathfrak{a}}\varepsilon^{\theta}, \; \;
        a_{1,i} \lvert z_{1,i}-z_{2,\pi(i)} \rvert
        \leq
        C_{\mathfrak{z}}\varepsilon^{\theta}, \; \textrm{for } i\in \mathfrak N,
    \end{align*}
    \begin{align*}
        a_{1, i}
        \leq
        C_{\mathfrak{a}}\varepsilon^{\theta}, \quad
        a_{2, j}
        \leq
        C_{\mathfrak{a}}\varepsilon^{\theta} \; \textrm{for } i\in \mathfrak N_{1}, \; j\in \mathfrak N_{2},
    \end{align*}
    where $\theta \in [\frac{1}{C}\eta^{2},C\eta^{2}]$ for some constant $C$ depending only on $R_{0}$. The constants $\varepsilon_{0}, C_{\mathfrak{a}} $ and $C_{\mathfrak{z}}$ depend only on $\Omega,q, k, \mathfrak p,\eta$.
\end{thm}

The stability of reconstructing the point sources is of H\"{o}lder type. Thus the inversion may be well-posed if the H\"{o}lder power $\theta$ is large enough. In fact, given explicitly in \eqref{theta} of the proof in Section \ref{section: stability for multiple point sources} the exponent $\theta$ depends on the minimum distance between the point sources $8\eta$ and tends to zero when $\eta$ goes to zero.

If the amplitudes of the point sources are uniformly strong, that is, there exists a strictly positive constant $\underline{a}$, such that $a_{1, i}> \underline{a}, i=1, \cdots, N_{1}$, one can easily recover a direct stability estimate on the location of the sources: $\lvert z_{1,i}-z_{2,\pi(i)} \rvert \leq \frac{C_{\mathfrak{z}}}{\underline{a}}\varepsilon^{\theta}, $ for $i \in \mathfrak N.$ We can also show the existence of a constant $\varepsilon_{1}>0$ such that $N_{1}=N_{2}$, and $\mathfrak N=\{1,\cdots, N_{1}\}, $ if $ 0\leq \varepsilon <\varepsilon_{1}$. Notice that in all the existing stability estimates for point sources the amplitudes were assumed to be uniformly strong and $ \varepsilon <\varepsilon_{1}$ (see for instance \cite{elbadia2012holder, elbadia2011inverse, ren2019imaging}, and references therein).

This paper is organized as follows. In Section \ref{section: forward problem Helmholtz with R.C.}, we study the existence and uniqueness of solutions to the forward problem. We also derive uniform bounds for the solutions when the source belongs to the set $\mathcal S$ in terms of the parameters $\mathfrak p$, $\Omega$, and $q$. Section \ref{section: stability for one point source} is devoted to the analysis of the inverse problem. Based on the CGO solutions, we first derive a H\"{o}lder type stability estimate for the inversion in the case where we have only one point source, that is $N=1$. Considering the geometrical assumption on the point sources we then use Carleman inequalities to estimate the solution of a Cauchy problem associated with the Helmholtz operator outside small balls containing the point sources. We finally derive a H\"{o}lder type stability for multiple point sources in Section \ref{section: stability for multiple point sources}. In the Appendix, we prove the existence of CGO solutions and provide the required Carleman inequality.

\section{The forward problem}\label{section: forward problem Helmholtz with R.C.}

In this section, we study the forward problem that consists of determining the radiated field $u$, the solution to the Helmholtz equation \eqref{equation: forward problem} from the knowledge of the source and the medium function. We present our principal result related to the forward problem as follows.

\begin{thm}\label{thm: uniform estimate forward}
    Assume that $S = \sum_{i=1}^{N} a_{i} \delta_{z_{i}}\in \mathcal S$. Then
    the inhomogeneous Helmholtz equation \eqref{equation: forward problem} admits a unique solution $u$. Moreover there exists a constant $\mathfrak c = \mathfrak c(\Omega,q, k, \varrho)>0,$
    such that
    \begin{align} \label{uniform bound for u forward}
    \lVert u \rVert_{H^{2}(\Omega^{4\varrho})} +\lVert u \rVert_{L^{2}(\Omega)} \leq \mathfrak c \lVert u_{0} \rVert_{L^{2}(\Omega)},
    \end{align}
    where $u_{0}(x) = \sum_{i=1}^{N}a_{i}\Phi(x,z_{i})$ with $\Phi$ is the free space Green function of the Helmholtz equation, and $\Omega^{\varrho} = \{x\in \Omega: \textrm{dist}(x, z_{i})> \varrho, i=1,\ldots,N\}.$
\end{thm}


\begin{proof}
    To analyze the singularity of the solution at the point sources, we shall split it into two parts. The first part is a singular explicit function $u_{0}$ and the second part is a smooth function $w$ that lies in $H^{2}(\Omega)$. Recall the expression of the free space Green function of the Helmholtz equation
    $$
        \Phi(x,z)=\frac{1}{4\pi}\frac{e^{ik\lvert x-z \rvert}}{\lvert x-z \rvert}, \quad x\not=z,
    $$
    and set $$u_{0}(x) =\sum_{i=1}^{N}a_{i}\Phi(x,z_{i}).$$ Then $w= u-u_{0}$ satisfies the following system
    \begin{align}\label{equation: perturbation w}
        \begin{cases}
            \Delta w + k^{2}(1+q(x)) w = -k^{2} q u_{0} \quad \text{in} \quad \mathbb{R}^{3}, \\
            \lim\limits_{r\to +\infty}r(\frac{\partial u}{\partial r}-iku)=0,\quad r=\left\|x \right\|,
        \end{cases}
    \end{align}

    Multiplying the equation above by $\Phi$, and integrating by parts yield

    \begin{align}\label{equation:w=Tn u}
        w(x)=-k^{2}\int_{\mathbb{R}^{3}}\Phi(x,y)q(y)u(y)dy.
    \end{align}

    Define the volume potential operator $V_{q}: L^{2}(\Omega)\longrightarrow H^{2}(\Omega)$, by $$V_{q}\phi(x) =k^{2}\int_{\Omega}\Phi(x,y)q(y)\phi(y)dy,\; \forall \phi\in L^{2}(\Omega).$$

    It is well known that $V_{q}$ is a bounded operator, that is
    \begin{align} \label{estimat: V Cop bound H2 L2}
        \lVert V_{q}\phi\rVert_{H^{2}(\Omega)}
        \leq k^{2}C_{op} \lVert q\rVert_{C^{0}(\overline \Omega)} \lVert \phi \rVert_{L^{2}(\Omega)},
    \end{align}
    where $C_{op}>0$ depends only on the domain $\Omega$ \cite{bao2010error, colton1998inverse}. Thus \eqref{equation:w=Tn u} can be reformulated as
    \begin{align}\label{equation:Reisz Fredholm u u0}
        (I+V_{q})u=u_{0},
    \end{align}
    where $I$ is the identity operator on $L^{2}(\Omega)$. Since $H^{2}(\Omega)$ is compactly embedded in
    $L^{2}(\Omega)$, the operator $V_{q}$ is compact from $L^{2}(\Omega)$ to itself. Using the Riesz-Fredholm theory this integral equation \eqref{equation:Reisz Fredholm u u0} is uniquely solvable if the kernel of the operator $I+V_{q}$ on $L^{2}(\Omega)$ is trivial. This can be proved using the unique continuation property for the Helmholtz equation \cite{colton1998inverse}. Since $u_{0} \in L^{2}(\Omega)$, \eqref{equation:Reisz Fredholm u u0} has a unique solution which in turn implies the existence and uniqueness of solutions to the main equation \eqref{equation: forward problem}. Therefore, the following estimate holds
    \begin{align*}
    \lVert u \rVert_{L^{2}(\Omega)} \leq \lVert (I+V_{q})^{-1}\rVert_{\mathcal L\left(L^{2}(\Omega)\right)}\lVert u_{0}\rVert_{L^{2}(\Omega)}.
    \end{align*}
    On the other hand we deduce from \eqref{equation:w=Tn u}-\eqref{estimat: V Cop bound H2 L2}, the following estimate
    \begin{align} \label{Estimate of w}
    \lVert w \rVert_{H^{2}(\Omega)} \leq k^{2}C_{op} \lVert q\rVert_{C^{0}(\overline \Omega)}\lVert (I+V_{q})^{-1}\rVert_{\mathcal L\left(L^{2}(\Omega)\right)}\lVert u_{0}\rVert_{L^{2}(\Omega)}.
    \end{align}
    Now, let $\chi\in C^{\infty}(\mathbb {R}^{3})$ be a cut off function satisfying $\chi= 1 $ on $\Omega^{2\varrho}$ and $\chi=0$ on $\mathbb {R}^{3} \setminus \overline{\Omega^{\varrho}}$. Since $S\in \mathcal S$, a simple calculation yields
    \begin{align} \label{Estimate0}
    \Delta (\chi u_{0}) = (-k^{2}\chi+\Delta \chi) u_{0} +2\nabla \chi \cdot \nabla u_{0}, \quad \textrm{in }\mathbb {R}^{3}.
    \end{align}
    Multiplying the previous equation by $\chi u_{0}$ and integrating by parts give
    \begin{align} \label{Estimate1}
    \lVert \nabla u_{0}\rVert_{L^{2}(\Omega^{2\varrho})}\leq
    \lVert \nabla (\chi u_{0})\rVert_{L^{2}(\Omega^{\varrho})} \leq \lVert -k^{2}\chi^{2}+\chi\Delta \chi-\frac{1}{2}\Delta \chi^{2} \rVert_{C^{0}(\overline{\Omega^{\varrho}})} \lVert u_{0}\rVert_{L^{2}(\Omega^{\varrho})}
    \end{align}
    Similarly, considering $v_{0}=\partial_{x_{j}}u_{0}$ and cut off function $\chi^{\prime}$ satisfying $\chi^{\prime}= 1 $ on $\Omega^{4\varrho}$ and $\chi^{\prime}=0$ on $\mathbb {R}^{3} \setminus \overline{\Omega^{2\varrho}}$  we have
    \begin{align}\label{Estimate2}
    \lVert \nabla v_{0}\rVert_{L^{2}(\Omega^{4\varrho})}\leq
    \lVert \nabla (\chi^{\prime} v_{0})\rVert_{L^{2}(\Omega^{2\varrho})} \leq C^{\prime} \lVert v_{0}\rVert_{L^{2}(\Omega^{2\varrho})}\leq C^{\prime}C^{\prime\prime} \lVert u_{0}\rVert_{L^{2}(\Omega^{\varrho})}
    \end{align}
    where $C^{\prime},C^{\prime\prime}$ are constants depending on $\Omega,\mathfrak{\varrho},k$ .Combining the \eqref{Estimate0}, \eqref{Estimate1}, \eqref{Estimate2}, we obtain
    \begin{align*}
    \lVert u_{0} \rVert_{H^{2}(\Omega^{4\varrho})} \leq c \lVert u_{0}\rVert_{L^{2}(\Omega)},
    \end{align*}
    where $c = c(\Omega, \varrho, k)>0$. Since $u = w+u_{0}$, using the previous inequality
    and estimate \eqref{Estimate of w}, taking
    $$
    \mathfrak c = (k^{2}C_{op}\lVert q\rVert_{C^{0}(\overline \Omega)} + 1) \lVert(I+V_{q})^{-1} \rVert_{\mathcal L\left(L^{2}(\Omega)\right)} + c +1,
    $$
    the proof is complete.
\end{proof}
We next study the stability of the solution $u$ to the Helmholtz equation under simultaneous perturbations of the medium function $q$ the source $S$.
\begin{proposition} \label{proposition 1}
    Let $S_{1} = \sum_{i=1}^{N_{1}} a_{1,i} \delta_{z_{1,i}}$ and $S_{2}= \sum_{i=1}^{N_{2}} a_{2,i} \delta_{z_{2,i}}$ be two sources in $\mathcal S$. Let $u_{1},u_{2}$ be the solutions to the inhomogeneous Helmholtz equation \eqref{equation: forward problem} with respectively point sources $S_{1} $ and $S_{2}$, and medium functions $q_{1}$ and $q_{2}$. Let $\mathfrak c_{1} = \mathfrak c(\Omega,q_{1}, k,\varrho)>0,$ and $\mathfrak c_{2} = \mathfrak c(\Omega,q_{2},k,\varrho)>0,$ where $\mathfrak c$ is constant in Theorem \ref{thm: uniform estimate forward} corresponding to respectively the medium functions $q_{1}$ and $q_{2}$. Then the following inequality holds

    \begin{align} \label{uniform bound for u}
    \lVert u_{1}-u_{2} \rVert_{H^{2}(\Omega^{\varrho})} +\lVert u_{1}-u_{2} \rVert_{L^{2}(\Omega)} \leq 
    \mathfrak c_{1} \mathfrak c_{2}
    \lVert u_{0,2} \rVert_{L^{2}(\Omega)} \lVert q_{1} -q_{2}\rVert_{C^{0}(\overline \Omega)}+ \mathfrak c_{1} \lVert u_{0,1} -u_{0,2} \rVert_{L^{2}(\Omega)},
    \end{align}
    where $u_{0,i}(x) = \sum_{j=1}^{N_{j}}a_{i,j}\Phi(x,z_{i, j}), \, i=1, 2,$ with $\Phi$ is the free space Green function of the Helmholtz equation, and $\Omega^{\varrho} = \{x\in \Omega: \textrm{dist}(x, z_{i,j})> \varrho,j=1,\ldots,N_{i},i=1,2\}.$
\end{proposition}
\begin{proof}
Let $w_{j} = u_{j} - u_{0,j}$ where $u_{0,j}=\sum_{i=1}^{N_{i}}a_{i}\Phi (x,z_{j,i})$. Analogous to \eqref{equation:Reisz Fredholm u u0} we deduce that $u_{j}$ satisfies
\begin{align*}
(I+V_{q_{j}})u_{j}=u_{0, j},
\end{align*}
for $j=1, 2$.  Therefore $u= u_{1}-u_{2}$ verifies

\begin{align*}
(I+V_{q_{1}})u = V_{q_{2}-q_{1}} u_{2} + u_{0, 1}-u_{0, 2}.
\end{align*}
Regarding the regularity of the map $V_{q},$ and $u_{2}$, we have
\begin{align*}
\lVert V_{q_{2}-q_{1}} u_{2} + u_{0, 1}-u_{0, 2}\rVert_{L^{2}(\Omega)} \leq
k^{2}C_{op}\lVert q_{1} -q_{2}\rVert_{C^{0}(\overline \Omega)} \lVert u_{2} \rVert_{L^{2}(\Omega)}+
\lVert u_{0,1} -u_{0,2} \rVert_{L^{2}(\Omega)},\\\leq
k^{2}C_{op} \mathfrak c_{2} \lVert u_{0,2}\rVert_{L^{2}(\Omega)} \lVert q_{1} -q_{2}\rVert_{C^{0}(\overline \Omega)}+
\lVert u_{0,1} -u_{0,2} \rVert_{L^{2}(\Omega)}.
\end{align*}
Concerning
\begin{align*}
    w_{1}-w_{2}=V_{q_{1}}u_{1}-V_{q_{2}}u_{2}=V_{q_{1}}(u_{1}-u_{2})+V_{q_{1}-q_{2}}u_{2},
\end{align*}
the rest of the proof follows the same steps of the one of Theorem \ref{thm: uniform estimate forward}.

\end{proof}


\begin{proposition}\label{lemma:regularity of f g depending on az} Fix $\alpha\in (0, \frac{1}{2})$.
    Let $S_{i} = \sum_{j=1}^{N} a_{i,j} \delta_{z_{i,j}}, i=1, 2,$ be two sources in $\mathcal S$, and let $u_{0,i}(x) = \sum_{j=1}^{N}a_{i,j}\Phi(x,z_{i, j}), \, i=1, 2,$ with $\Phi$ is the free space Green function of the Helmholtz equation. Let $\pi^{*}$ be a permutation acting on $\{1,2,\ldots,N\}$. Then,
    the following inequality holds
    \begin{align*}
    \lVert u_{0, 1}-u_{0, 2} \rVert_{L^{2}(\Omega)} \leq \mathfrak m_{1} \sum_{i=1}^{N} |a_{1,i} - a_{2,\pi^{*}(i)}|
    + \mathfrak m_{2} \sum_{i=1}^{N} |z_{1,i} - z_{2,\pi^{*}(i)}|^\alpha,
    \end{align*}
    with $$ \mathfrak m_{1} = \max_{z\in \overline \Omega} \lVert \Phi(\cdot, z)
        \rVert_{L^{2}(\Omega)}, \; 
        \mathfrak m_{2}= 4\pi^{\frac{1}{2}}(1+(2R_{0})^{\frac{3}{4}})^{\frac{1}{2}}
        + 2\pi^{\frac{1}{2}} \lvert z_{1,i}-z_{2,\pi^{*}(i)} \rvert^{1-\frac{3\alpha}{2}}. $$
\end{proposition}

\begin{proof} By the triangle inequality, we have
    $$
        \lVert u_{0, 1}-u_{0, 2} \rVert_{L^{2}(\Omega)} \leq
        \sum_{i=1}^{N} |a_{1,i} - a_{2,\pi^{*}(i)} |\lVert \Phi(\cdot, z_{1,i})
        \rVert_{L^{2}(\Omega)}+|a_{2,\pi^{*}(i)}| |\lVert \Phi(\cdot, z_{1,i})
        - \Phi(\cdot, z_{2,\pi^{*}(i)}) \rVert_{L^{2}(\Omega)}.
    $$
    On the other hand, dividing the integral domain $\Omega$ of $\Phi(\cdot, z_{1,i})
    - \Phi(\cdot, z_{2,\pi^*(i)})$ into two parts, $\Omega_{S}$ and $\Omega\backslash\Omega_{S}$, 
    where
    $$\Omega_{S}=\bigl\{x\in \Omega \big| \lvert x-z_{1,i}\rvert\leq \lvert z_{1,i}-z_{2,\pi^{*}(i)}\rvert^{\alpha} \bigr\} \cup \bigl\{ x\in \Omega \big| \lvert x-z_{2,\pi^{*}(i)}\rvert\leq \lvert z_{1,i}-z_{2,\pi^{*}(i)}\rvert^{\alpha}\bigr\},$$ 
    a forward computation yields
    \begin{align*}
        &\lVert \Phi(\cdot, z_{1,i})
        - \Phi(\cdot, z_{2,\pi^*(i)}) \rVert_{L^{2}(\Omega)} \\
        \leq
        & 
        4\pi^{\frac{1}{2}}(1+(2R_{0})^{\frac{3}{4}})^{\frac{1}{2}}\lvert z_{1,i}-z_{2,\pi^{*}(i)}\rvert^{\frac{\alpha}{2}}
        + 2\pi^{\frac{1}{2}} \lvert z_{1,i}-z_{2,\pi^{*}(i)} \rvert^{1-\alpha}.
    \end{align*}
    Combining the above inequalities provides the desired estimate.
\end{proof}
\section{The inverse problem}\label{section: inverse problem}
The inverse problem is to determine the point sources $S $ from the knowledge of the radiated field $ u|_{\partial\Omega}$ assuming that the medium function $q(x)$ is given. Our goal is to study the stability of the inverse problem, that is, to estimate the closeness of the two sources $S_{1}$ and $S_{2}$ in some suitable metric if the difference between their radiated fields $u_{1}$ and $u_{2}$ are close on the boundary of $\Omega$. We first follow the approach of \cite{bao2021recovering} based on CGO solutions to derive a stability estimate for the case of a single point source. Then, using Carleman inequalities for the Helmholtz operator, we derive a stability estimate to the Cauchy problem in the domain excluding small balls with a fixed radius containing the point sources. Afterwards, we show that if $ u_{1}|_{\partial\Omega}$ and $ u_{2}|_{\partial\Omega}$ are close enough, each ball contains exactly two point sources. We finally derive a H\"{o}lder type stability for multiple point sources by applying the previous results to each of the small balls respectively.

\subsection{Stability for the inverse problem with one point source}\label{section: stability for one point source} We first consider one point source and then extend the result to multiple point sources.

\begin{thm}\label{thm: stability a z}
    Let $S_{1} = a_{1} \delta_{z_{1}}$ and $S_{2}=
        a_{2} \delta_{z_{2}}$ be two sources in $\mathcal S$.
    Let $u_{1},u_{2}$ be the solutions to the inhomogeneous Helmholtz equation \eqref{equation: forward problem} with respectively point sources $S_{1} $ and $S_{2}$, and the same medium function $q$. Denote
    $$
        \epsilon=\lVert u_{1}-u_{2}\rVert_{H^{\frac{1}{2}}(\partial\Omega)} + \lVert \frac{\partial u_{1}}{\partial\nu}- \frac{\partial u_{2}}{\partial\nu} \rVert_{H^{-\frac{1}{2}}(\partial\Omega)}
    $$
    Then the following stability estimate
    \begin{align*}
        \lvert a_{1}-a_{2} \rvert\leq C_{a}\epsilon, \;\;
        a_{1} \lvert z_{1}-z_{2} \rvert\leq C_{z}\epsilon,
    \end{align*}
    holds, where
    \begin{align*}
         & C_{a}=2e^{7\sqrt{M^{2}+4M_{G}^{2}C_{L}^{2}}R_{0}}\lvert \partial\Omega\rvert^{\frac{1}{2}},\;\; C_{z}= M_{G}^{-1}C_{L}^{-1}4e^{11\sqrt{M^{2}+4M_{G}^{2}C_{L}^{2}}R_{0}} \lvert \partial\Omega\rvert^{\frac{1}{2}}, \\
         & M=\max\{\frac{2916R_{0}k^{2}}{\pi}\lVert q \rVert_{H^{2}(\Omega)}
        C_{s} ,\frac{2R_{0}k^{2}}{\pi}\lVert q \rVert_{C^{2}(\Omega)}, \frac{24}{11}M_{G}C_{L},\frac{\pi}{R_{0}},k \},                                                                     \\
         & M_{G}=\left(\left(\frac{R_{0}k^{2}}{\pi}(4\lVert q\rVert_{C^{0}(\overline D)}+24)+\frac{27}{2}\right)
         \frac{18R_{0}k^{2}}{\pi |\Im\xi|} + \frac{4 R_{0}k^{2}}{\pi} \right)\lVert q \rVert_{H^{2}(D)}
    \end{align*}
    and $C_{L}>0$ is the Sobolev embedding constant as stated in Remark \ref{remark:1/12 gradient decay and continous dependence on xi}.
\end{thm}

\begin{proof}

    It is easy to check that $U=u_{1}-u_{2}$ satisfies the following system
    \begin{align*}
        \begin{cases}
            \Delta U + k^{2}(1+q(x)) U = a_{1} \delta_{z_{1}} -a_{2} \delta_{z_{2}} \quad \textrm{in } \mathbb {R}^{3}, \\
            \lim\limits_{r\to +\infty}r(\frac{\partial U}{\partial r}-ikU)=0,\quad r=\left\|x \right\|,
        \end{cases}
    \end{align*}

    Construct a CGO solution $v_{\xi}(x)$ satisfying
    \begin{align*}
        \Delta v_{\xi} + k^{2}(1+q(x))v_{\xi} = 0,
    \end{align*}
    as in Lemma \ref{lemma:cgo decay n2} where $v_{\xi}$ has the form $v_{\xi}(x)=e^{ix\cdot\xi}(1+\phi_{\xi}),\xi\in \mathbb{C}^{3},\xi\cdot\xi=k^{2}$. Here we choose $\xi$ with a specific form
    $$
        \xi=i(t_{1}e_{1}+t_{2}e_{2})+\sqrt{k^{2}+{t_{1}}^{2}+{t_{2}}^{2}} e_{3},
    $$
    where $\{e_{1},e_{2},e_{3}\}$ is a orthonormal Cartesian basis of $\mathbb{R}^{3}$, and $t_{1},t_{2}\in \mathbb{R}$.

    Use the second Green identity for $U$ and $v_{\xi}$ we have
    \begin{align}\label{equation:gap reciprocity identity}
        a_{1}v_{\xi}(z_{1})-a_{2}v_{\xi}(z_{2})=\int_{\partial\Omega}\left(\frac{\partial u_{1}}{\partial \nu}- \frac{\partial u_{2}}{\partial \nu}\right)v_{\xi}
        - (u_{1}-u_{2})\frac{\partial v_{\xi}}{\partial \nu} dS.
    \end{align}

    The stability estimates of $a_{1}-a_{2}$ and $z_{1}-z_{2}$ are based on a specific choice of the CGO parameters $t_{1},t_{2},e_{1},e_{2},e_{3}$ and via direct calculations in three steps.

    \paragraph{Step 1}
    Control the CGO potential $\phi_{\xi}$ in $C^{0}(\overline \Omega)$. Choose $t_{1}=t^{*}_{1}=M$, then $\lvert \Im \xi \rvert\geq t_{1}\geq M$, Remark \ref{remark:1/12 gradient decay and continous dependence on xi} implies for any $t_{2} \in \mathbb{R}$ and orthogonal triplet $e_{1},e_{2},e_{3}$ the corresponding CGO potential $\phi_{\xi}$ satisfies $\lVert\phi_{\xi} \rVert_{C^{0}(\overline \Omega)} \leq \frac{1}{12},\lVert \nabla \phi_{\xi} \rVert_{C^{0}(\overline \Omega)} \leq M_{G}C_{L}$.

    \paragraph{Step 2}
    Estimate $a_{1}-a_{2}$. Given $t^{*}_{1}=M$, we would like to find $t_{2}^{*}\in \mathbb{R}$ and $e_{1}^{*},e_{2}^{*},e_{3}^{*}$ such that for $$\xi^{*}=i(t_{1}^{*}e_{1}^{*}+t_{2}^{*}e_{2}^{*})+\sqrt{k^{2}+{t_{1}^{*}}^{2}+{t_{2}^{*}}^{2}} e_{3}^{*},$$
    we have $\lvert v_{\xi^{*}}(z_{1})\rvert  = \lvert v_{\xi^{*}}(z_{2}) \rvert $. This latter associated with \eqref{equation:gap reciprocity identity} and Step 1 and $\lvert \xi\rvert,\lvert\nabla \phi_{\xi} \rvert \leq \sqrt{3}\lvert \Im\xi\rvert$ implies

    \begin{align}\label{estimate: a1-a2 bounded by 1/v boundary}
    \lvert a_{1} -a_{2}\rvert = \frac{1}{\lvert v_{\xi^{*}}(z_{1})\rvert} \bigl\lvert a_{1} \lvert v_{\xi^{*}}(z_{1}) \rvert  -a_{2} \lvert v_{\xi^{*}}(z_{2})\rvert \bigr\rvert \leq  \frac{1}{\lvert v_{\xi^{*}}(z_{1})\rvert} \left \lvert a_{1}  v_{\xi^{*}}(z_{1})   -a_{2}  v_{\xi^{*}}(z_{2}) \right\rvert \nonumber
    \\
    \leq 
    \frac{1}{\lvert v_{\xi^{*}}(z_{1})\rvert} \left( \lVert u_{1}-u_{2}\rVert_{H^{\frac{1}{2}}(\partial\Omega)} + \left\lVert
    \frac{\partial u_{1}}{\partial \nu}- \frac{\partial u_{2}}{\partial \nu} \right\rVert_{H^{-\frac{1}{2}}(\partial\Omega)} \right) e^{5R_{0}|\Im\xi^{*}|}\lvert\partial\Omega\rvert^{1/2},
    \end{align}

    which is close to the desired result. Next, we show the existence of $\xi^{*}$ with such properties.

    If $z_{1}=z_{2}$,  $\lvert v_{\xi^{*}}(z_{1})\rvert  = \lvert v_{\xi^{*}}(z_{2}) \rvert $ is trivial. So we only have to handle the case where $z_{1}\neq z_{2}$, and find $t_{2}^{*},e_{1}^{*},e_{2}^{*},e_{3}^{*}$ such that $ \left\lvert \frac{v_{\xi^{*}}(z_{1})}{v_{\xi^{*}}(z_{2})} \right\rvert =1.$

    Firstly, we can choose $e_{0}^{*}$ as an unit vector satisfying $e_{0}^{*}\times (z_{1}-z_{2})\neq 0$ and an orthogonal triplet $\{ e_{1}^{*},e_{2}^{*},e_{3}^{*} \} \subset \mathbb{R}^{3}$ is constructed as following  
    $$
        e_{2}^{*}=\frac{z_{1}-z_{2}}{\lvert z_{1}-z_{2} \rvert},
        e_{3}^{*}=e_{0}^{*}\times e_{2}^{*},
        e_{1}^{*}=e_{3}^{*}\times e_{2}^{*}.
    $$
    Note that $\xi=i(t_{1}^{*}e_{1}^{*}+t_{2}e_{2})^{*}+\sqrt{k^{2}+{t_{1}^{*}}^{2}+{t_{2}}^{2}} e_{3}^{*}$ is a function of $t_{2}$ and $(z_{1}-z_{2})\cdot e_{1}^{*}=0,(z_{1}-z_{2})\cdot e_{3}^{*}=0.$

    Secondly, we introduce an auxiliary function of $t_{2}$ as follows.
    $$
    F(t_{2})=a(t_{2})e^{-t_{2}\lvert z_{1}-z_{2}\rvert},
    $$
    where 
    $$
    a(t_{2})=\left\lvert \frac{1+\phi_{\xi}(z_{1})}{1+\phi_{\xi}(z_{2})}\right\rvert \in [\frac{1-1/12}{1+1/12},\frac{1+1/12}{1-1/12}]=[\frac{11}{12},\frac{13}{12}].
    $$ 
    So $F(t_{2})\to\pm \infty$ as $t_{2}\to \mp \infty$. Notice that $\left\lvert \frac{v_{\xi^{*}}(z_{1})}{v_{\xi^{*}}(z_{2})} \right\rvert = F(t_{2})$, the mean-value Theorem implies the existence of a solution of $F(t_{2}^{*})=1$.

    Moreover, with $\zeta=\frac{\phi_{\xi}(z_{1})-\phi_{\xi}(z_{2})}{1+\phi_{\xi}(z_{2})}$ and $\lvert \zeta\rvert\leq 2/11\leq 1$ we have $\lvert \ln (a(t_{2}))\rvert \leq \lvert \ln (\lvert 1+ \lvert\zeta \rvert)\rvert \leq  \lvert \zeta \rvert $, then 
    $$
    \lvert \ln(a(t_{2}))\rvert \leq  \frac{\lvert\phi_{\xi}(z_{1})-\phi_{\xi}(z_{2}) \rvert }{ 1 + \lvert \phi_{\xi}(z_{2})\rvert}\leq 2\lvert\phi_{\xi}(z_{1})-\phi_{\xi}(z_{2}) \rvert.
    $$

    So $F(t_{2}^{*})=1$, which is equivalent to $t_{2}^{*}=\frac{1}{\lvert z_{1}-z_{2} \rvert}\ln(a(t_{2}^{*}))$, implies
    $$
    \lvert t_{2}^{*} \rvert \leq \frac{\lvert \ln(a(t_{2}))\rvert}{\lvert z_{1}-z_{2} \rvert}  
    \leq \frac{ 2\lvert \phi_{\xi}(z_{1})-\phi_{\xi}(z_{2}) \rvert}{\lvert z_{1}-z_{2} \rvert} \leq 2M_{G}C_{L}.
    $$
    Hence, we can find $t_{2}^{*} \in [-2M_{G}C_{L},2M_{G}C_{L}]$ such that $|v_{\xi^{*}}(z_{1})|= |v_{\xi^{*}}(z_{2})|$ where
    $$
        \xi^{*}=i(t_{1}^{*}e_{1}^{*}+t_{2}^{*}e_{2}^{*}) + \sqrt{k^{2}+{t_{1}^{*}}^{2}+{t_{2}^{*}}^{2}} e_{3}^{*},
        t_{1}^{*}=M,\lvert\Im\xi\rvert\leq \sqrt{M^{2}+4M_{G}^{2}C_{L}^{2}}
    $$


    From the identity \eqref{estimate: a1-a2 bounded by 1/v boundary} and choosing $\xi^{*}$ as above, we have
    \begin{align}\label{estimate:a1-a2 epsilon}
        \lvert a_{1}-a_{2} \rvert
         & \leq
        \left( \lVert u_{1}-u_{2}\rVert_{H^{\frac{1}{2}}(\partial\Omega)} + \left\lVert
        \frac{\partial u_{1}}{\partial \nu}- \frac{\partial u_{2}}{\partial \nu}\right\rVert_{H^{-\frac{1}{2}}(\partial\Omega)} \right) 
        2e^{7\sqrt{M^{2}+4M_{G}^{2}C_{L}^{2}}R_{0}}\lvert \partial\Omega\rvert^{\frac{1}{2}}\nonumber
        \\ &\leq 
        2e^{7\sqrt{M^{2}+4M_{G}^{2}C_{L}^{2}}R_{0}}\lvert \partial\Omega\rvert^{\frac{1}{2}}\epsilon
    \end{align}

    \paragraph{Step 3}
    Estimate $z_{1}-z_{2}$. We recall the following identity
    \begin{align*}
        a_{1}v_{\xi}(z_{1})-a_{1}v_{\xi}(z_{2})+a_{1}v_{\xi}(z_{2})-a_{2}v_{\xi}(z_{2})
        =\int_{\partial\Omega}\left(\frac{\partial u_{1}}{\partial \nu}- \frac{\partial u_{2}}{\partial \nu}\right)v_{\xi}
        - (u_{1}-u_{2})\frac{\partial v_{\xi}}{\partial \nu} dS,
    \end{align*}
    and a direct consequence is that

    \begin{align} \label{equation: z1-z2 gradient calculation}
        a_{1}\lvert v_{\xi}(z_{1})-v_{\xi}(z_{2}) \rvert
        \leq
        \left \lvert \int_{\partial\Omega}\left(\frac{\partial u_{1}}{\partial \nu}- \frac{\partial u_{2}}{\partial \nu}\right)v_{\xi}
        - (u_{1}-u_{2})\frac{\partial v_{\xi}}{\partial \nu} dS \right \rvert
        + \lvert a_{1}-a_{2} \rvert \lvert v_{\xi}(z_{2})\rvert.
    \end{align}

    If we can derive a lower bound of $|v_{\xi}(z_{1})-v_{\xi}(z_{2})|$, then combined with \eqref{estimate:a1-a2 epsilon}, we obtain the desired estimate. We now provide details of this strategy with a specific choice of $t_{2},e_{1},e_{2},e_{3}$.

    If $z_{1}=z_{2}$, the stability estimate in Theorem \ref{thm: stability a z} becomes trivial, so we can further consider the case $z_{1}\neq z_{2}$. Choose the CGO parameter for $\xi^{*}$ as
    \begin{align*}
        t_{2}^{*}=t_{1}^{*}=M,
        e_{2}^{*}=\frac{z_{1}-z_{2}}{\lvert z_{1}-z_{2} \rvert},
        e_{3}^{*}=e_{0}^{*}\times e_{2}^{*},
        e_{1}^{*}=e_{3}^{*}\times e_{2}^{*}.
    \end{align*}
    where $e_{0}^{*}$ is an unit vector satisfying $e_{0}\times (z_{1}-z_{2})\neq 0.$ as we did in Step 2.
    Constructing an auxiliary function of $s$ as below,
    $$
        H(s)=v(z_{2}+s(z_{1}-z_{2}))=e^{i(z_{2}+s(z_{1}-z_{2}))\cdot\xi^{*}}(1+\phi_{\xi^{*}}(z_{2}+s(z_{1}-z_{2}))),
    $$
    A direct calculation gives
    \begin{align*}
        H^{\prime}(s)
         & =i\xi^{*} e^{i(z_{2}+s(z_{1}-z_{2}))\cdot\xi^{*}} (z_{1}-z_{2})\cdot \Bigl(1+\phi_{\xi^{*}}(z_{2}+s(z_{1}-z_{2}))\Bigr)
        +e^{i(z_{2}+s(z_{1}-z_{2}))\cdot \xi^{*}}\nabla\phi_{\xi^{*}}\cdot(z_{1}-z_{2}),                                                 \\
         & = e^{i(z_{2}+s(z_{1}-z_{2}))\cdot\xi^{*}} \Biggl( -M\lvert z_{1}-z_{2} \rvert \biggl(1+\phi_{\xi^{*}}\bigl(z_{2}+s(z_{1}-z_{2})\bigr)\biggr)+ \nabla\phi_{\xi^{*}}\cdot(z_{1}-z_{2})\Biggr).
    \end{align*}

    Notice that we have $\lVert \phi_{\xi} \rVert_{C^{0}(\overline \Omega)} \leq \frac{1}{12},\lVert\nabla \phi_{\xi} \rVert_{C^{0}(\overline \Omega)} \leq M_{G}C_{L}$, and $ \frac{24}{11}M_{G}C_{L}\leq = M $, the mean value theorem for $\Re H^{\prime}$ implies
    $$
        \lvert \Re( v_{\xi}(z_{1}) - v_{\xi}(z_{2})) \rvert \geq \min_{s\in (0,1)}\lvert \Re H^{\prime}(s) \rvert
        \geq (\frac{11}{12} M \lvert z_{1}-z_{2}\rvert
        -M_{G}C_{L}\lvert z_{1}-z_{2}\rvert ) e^{-2t_{2}^{*}R_{0}}
        \geq M_{G}C_{L}\lvert z_{1}-z_{2}\rvert e^{-2t_{2}^{*}R_{0}}.
    $$

    Combining the inequality above with \eqref{estimate:a1-a2 epsilon}, and \eqref{equation: z1-z2 gradient calculation}, we obtain
    \begin{align*}
        |a_{1}| \lvert z_{1}-z_{2}\rvert
         & \leq
        M_{G}^{-1}C_{L}^{-1}e^{2t_{2}^{*}R_{0}}
        \left(
        \left \lvert \int_{\partial\Omega}\left(\frac{\partial u_{1}}{\partial \nu}- \frac{\partial u_{2}}{\partial \nu}\right)v_{\xi}
        - (u_{1}-u_{2})\frac{\partial v_{\xi}}{\partial \nu} dS \right \rvert
        + \lvert a_{1}-a_{2} \rvert \lvert v_{\xi}(z_{2})\rvert
        \right),
        \\
         & \leq
        M_{G}^{-1}C_{L}^{-1}e^{2t_{2}^{*}R_{0}}
        \left(
        2e^{5t_{2}^{*}R_{0}}\lvert \partial\Omega\rvert^{\frac{1}{2}}\epsilon
        + \lvert a_{1}-a_{2} \rvert e^{2t_{2}^{*}R_{0}}
        \right),                    \\
         & \leq M_{G}^{-1}C_{L}^{-1}
        4e^{11\sqrt{M^{2}+4M_{G}^{2}C_{L}^{2}}R_{0}} \lvert \partial\Omega\rvert^{\frac{1}{2}}\epsilon.
    \end{align*}
\end{proof}

\section{Stability for the inverse problem with multiple point sources: Proof of Theorem \ref{thm: stability for multiple a z}}\label{section: stability for multiple point sources}

We begin with three auxiliary results that are of interest themselves. The first result is related to the distribution of sources in the domain $\Omega$.

\begin{lemma} \label{lemma Geometry}
    Let $S_{1} = \sum_{i=1}^{N_{1}} a_{1,i} \delta_{z_{1,i}}$ and $S_{2}= \sum_{i=1}^{N_{2}} a_{2,i} \delta_{z_{2,i}}$ be two sources in $\mathcal S$. Then the ball $B_{3\eta}(z_{i,j})$ contains a maximum of two
    point sources.
\end{lemma}
\begin{proof}
The result is a direct consequence of the assumption $\lvert z_{i,j}-z_{i,k} \rvert \geq 8\eta$ for $i=1, 2,$ and
$ 1\leq j\neq k \leq N_{i}.$
\end{proof}

Recall the definition $\mathfrak N \subset \{1, \cdots, N_{1}\}$ is the set for which $B_{3\eta}(z_{1,j}), j \in \mathfrak N,$ contains exactly two point sources, and $\pi: \mathfrak N \rightarrow \{1, \cdots, N_{2}\}$ is the one-to-one application such that $z_{2,\pi(j)} \in B_{3\eta}(z_{1,j})$ for all $j\in \mathfrak N,$ and $\mathfrak N_{1} = \{1, \cdots, N_{1}\} \setminus \mathfrak N,\, \mathfrak N_{2} = \{1, \cdots, N_{2}\} \setminus \pi(\mathfrak N). $

Denote
\begin{align}\label{domain Beta1}
\mathfrak B_{\eta, 1} = \left(\cup_{i\in \mathfrak N}B_{\eta}(z_{1,i})\right) \cup \left(\cup_{i= \mathfrak N_{1}}B_{\eta}(z_{1,i})\right)\cup\left(\cup_{i\in \mathfrak N_{2}}B_{\eta}(z_{2,i})\right), \\
\label{domain Beta2}
\mathfrak B_{\eta, 2} = \left(\cup_{i\in \mathfrak N}B_{2\eta}(z_{1,i})\right) \cup \left(\cup_{i= \mathfrak N_{1}}B_{2\eta}(z_{1,i})\right)\cup\left(\cup_{i\in \mathfrak N_{2}}B_{2\eta}(z_{2,i})\right),\\
\mathfrak B_{\eta, 3} = \left(\cup_{i\in \mathfrak N}B_{3\eta}(z_{1,i})\right) \cup \left(\cup_{i= \mathfrak N_{1}}B_{3\eta}(z_{1,i})\right)\cup\left(\cup_{i\in \mathfrak N_{2}}B_{3\eta}(z_{2,i})\right).
\end{align}

By the construction above, $\mathfrak B_{\eta, 1}$ and $ \mathfrak B_{\eta, 2}$ are sets of disjoint balls. Therefore the open sets
\begin{align} \label{domain eta }
\Omega_{\eta, i} = \Omega \setminus \overline{\mathfrak B_{\eta, i}}, \; i=1,2, 3,
\end{align}
are smooth and connected, and we have $$\Omega_{\eta, 3} \subsetneqq
    \Omega_{\eta, 2} \subsetneqq \Omega_{\eta, 1}.$$


\begin{lemma}\label{lemma:cut-off funciton}
    There exists a cut-off function $\chi\in C^{\infty}(\Omega),0\leq \chi \leq 1,$ such that
    \begin{align*}
        \begin{cases}
            \chi\equiv 0 \quad in \quad \mathfrak B_{\eta, 1}, \\
            \chi\equiv 1 \quad in \quad \Omega \setminus \overline{\mathfrak B_{\eta, 2}},
        \end{cases}
    \end{align*}
    with $\lVert \chi \rVert_{C^{2}(\overline{\Omega})} \leq C\frac{\left(N_{1}+N_{2}\right)^{2}}{\eta^{2}}$ where $C>0$ is an universal constant independent of $\Omega$.
\end{lemma}
\begin{proof}
    Let $z^{\star} \in \Omega$ satisfying $dist(z^{\star}, \partial\Omega)\geq 8\eta$. We next show the existence of a cut off function $\chi_{z^{\star}}\in C^{\infty}(\Omega)$ verifying $\lvert \chi_{z^{\star}}\rvert \leq 1,\lvert \nabla\chi_{z^{\star}}\rvert\leq \frac{C}{\eta},\lvert \bigtriangleup\chi_{z^{\star}}\rvert\leq \frac{C}{\eta^{2}}$, and
    \begin{align*}
        \begin{cases}
            \chi_{z^{\star}} \equiv 0 \quad in \quad B_{\eta}(z^{\star}) \\
            \chi_{z^{\star}} \equiv 1 \quad in \quad \mathbb{R}^{3}\setminus\overline{B_{2\eta}(z^{\star})}.
        \end{cases}
    \end{align*}

    The construction of $\chi_{z^{\star}}$ is well-known, for completeness an explicit expression for $\chi_{z^{\star}}$ is shown in two steps.

    \paragraph{Step 1}
    We first introduce a cut-off function $f$, and a smooth mollifier $g$ as follows
    \begin{align*}
        \begin{cases}
            f \equiv 0 \quad \textrm{in } \quad B_{\frac{3\eta}{2}}(z^{\star}), \\
            f \equiv 1 \quad \textrm{in } \quad \mathbb{R}^{3} \backslash\overline{B_{\frac{3\eta}{2}}(z^{\star})},
        \end{cases}
        \begin{cases}
            g = \left(\int_{B_{1}(O)} e^{\frac{1}{\lvert s \rvert^{2}-1}} ds\right)^{-1}
            e^{\frac{1}{\lvert x \rvert^{2}-1}} \quad \textrm{in } B_{1}(O), \\
            g \equiv 0 \quad \textrm{in } \quad \mathbb{R}^{3} \backslash \overline{B_{1}(O)}.
        \end{cases}
    \end{align*}
    The function $g$ lies in $C_{0}^{\infty}(\mathbb {R}^{3})$ and satisfies $\int_{\mathbb {R}^{3}} g dx =1$. For $0\leq a \leq 1 $ let $g_{a}(x)=g(\frac{x}{a})\frac{1}{a^{3}}$. Then $\int_{\mathbb{R}^{3}} g_{a} dx =1 , \lVert g_{a}\rVert_{C^{2}(\overline{B_{1}(O)})} \leq \frac{1}{a^{2}}\lVert g \rVert_{C^{2}(\overline{B_{1}(O)})}\leq \frac{m_{g}}{a^{2}}$ where $m_{g}>0$ is a constant that depends only on $g$.

    \paragraph{Step 2}
    Let $\chi_{z^{\star}}=f\ast g_{\frac{\eta}{2}}$, where $\ast$ stands for the convolution operation between two functions on $\mathbb {R}^{3}$. It is straightforward that $\chi_{z^{\star}}\in C^{\infty}(\Omega)$. Next, we check that $\chi_{z^{\star}}$ has the required properties.

    If $x\in B_{\eta}(z^{\star})$ then $\chi_{z^{\star}}(x)=\int_{0\leq \lvert y-x \rvert\leq \frac{\eta}{2} }f(y)g_{\frac{\eta}{2}}(y-x)dy=0$. While if $x\in \Omega\backslash\overline{B_{2\eta}(z^{\star})}$ then $\chi_{z^{\star}}(x)=\int_{0\leq \lvert y-x \rvert\leq \frac{\eta}{2} }f(y)g_{\frac{\eta}{2}}(y-x)dy=1$.

    Note that $\partial_{i}\chi_{z^{\star}}=f*\partial_{i}g_{\frac{\eta}{2}}=\frac{2}{\eta}f*\partial_{i}g,\, \partial_{i}\partial_{j}\chi_{z^{\star}}=f*\partial_{i}\partial_{j}g_{\frac{\eta}{2}}=\frac{4}{\eta^{2}}f*\partial_{i}\partial_{i}g$. Then the construction of $\chi_{z^{\star}}$ is completed.

    Consider $\chi(x)=\prod_{i\in \mathfrak N}\chi_{z_{1,i}}(x) \prod_{i\in \mathfrak N_{1}} \chi_{z_{1,i} }(x)\prod_{i\in \mathfrak N_{2}} \chi_{z_{2,i}}(x), \; x\in \Omega, $ and a direct calculation gives the desired
    result.
\end{proof}

\begin{lemma}\label{lemma:Carleman Estimate}
    Under the same assumption of Theorem \ref{thm: stability for multiple a z}, let $u = u_{1} -u_{2}$. $\varepsilon = \lVert u \rVert_{H^{1}(\partial\Omega)} + \lVert \frac{\partial u}{\partial \nu} \rVert_{L^{2}(\partial\Omega)}.$ There exists a constant $\varepsilon_{0}>0,\tau_{1}>0,M_{u}$, that only depends on $\Omega, q, k,$ and $ \mathfrak p$, such that

    If $\varepsilon < \varepsilon_{0},$ then
    \begin{align*}
    \lVert u \rVert_{H^{1}(\Omega_{\eta, 3})}\leq \sqrt{2} \left(\frac{M_{u}}{\eta}\right)^{1-{\tilde \theta}} \varepsilon^{\tilde \theta},
    \end{align*}
    with ${\tilde \theta}= \frac{5\eta^{2}}{2+R_{0}^{2}-4\eta^{2}}$.

    If $\varepsilon \geq \varepsilon_{0},$ then
    \begin{align*}
    \lVert u \rVert_{H^{1}(\Omega_{\eta, 3})}
    \leq \sqrt{2} C(\eta) \varepsilon.
    \end{align*}
    
\end{lemma}

\begin{proof}
    The result is a direct consequence of Carleman estimates for the Helmholtz equation \cite{hrycak2004increased, choulli2016applications}. A H\"{o}lder stability of the Cauchy problem has been proved in \cite{alessandrini2009stability, choulli2016applications}.

    We further assume that $1\in \mathfrak N$. The other case $1\in \mathfrak N_{1}$ can be treated similarly. Since $z_{1,1} \notin \Omega_{\eta, 1}$, $\phi=\lvert x-z_{1,1}\rvert^{2}$ has no critical points
    in $\Omega_{\eta, 1}$, and can be used as weight function in Proposition \ref{Carleman estimate} with $\omega = \Omega_{\eta, 1}$ and $Au=\Delta u + k^{2}(1+q)u=\sum_{i=1}^{N_{1}}a_{1,i}\delta_{1,i}-\sum_{i=1}^{N_{2}}a_{2,i}\delta_{z_{2,i}}$ satisfying $Au=0$ in $\Omega_{\eta,1}$. We have for $\tau\geq \tau_{0}\geq 2$
    \begin{align*}
         & \tau^{2}\lVert e^{\tau\phi} v \rVert^{2}_{L^{2}(\Omega_{\eta, 1})} + \tau\lVert e^{\tau\phi} \nabla v\rVert^{2}_{L^{2}(\Omega_{\eta, 1})}                                                                                   \\
         & \leq C\left( \lVert e^{\tau\phi}A v \rVert^{2}_{L^{2}(\Omega_{\eta, 1})} + \tau^{3}\lVert e^{\tau\phi}v \rVert^{2}_{L^{2}(\partial\Omega_{\eta, 1})} + \tau\lVert e^{\tau\phi}\nabla v \rVert^{2}_{L^{2}(\partial\Omega_{\eta, 1})} \right).
    \end{align*}

    Let $v=\chi u$ where $\chi$ is a cut-off function as stated in Lemma \ref{lemma:cut-off funciton} which satisfies $\chi=0$ in $\mathfrak{B_{\eta,1}}$ and $\chi=1$ in $\Omega\backslash\overline{\mathfrak{B_{\eta,2}}}$, then the boundary integral over $\mathfrak B_{\eta, 1}$ disappears. Expanding the term $A(\chi u)$ into four terms applying $\chi Au = 0$ in $\Omega_{\eta,1}$ and shrinking the integral domain for the left-hand side term above from $\Omega_{\eta, 1}$ to $\Omega_{\eta, 2}$, we have
    \begin{align*}
         & \tau^{2}\lVert e^{\tau\phi} u\rVert^{2}_{L^{2}(\Omega_{\eta, 2})} + \tau\lVert e^{\tau\phi} \nabla u\rVert^{2}_{L^{2}(\Omega_{\eta, 2})}                                                                                                                                                                        \\
         & \leq C\left( \lVert e^{\tau\phi}\nabla\chi\cdot\nabla u \rVert^{2}_{L^{2}(\Omega_{\eta, 1})} + \lVert e^{\tau\phi}u\Delta \chi \rVert^{2}_{L^{2}(\Omega_{\eta, 1})} + \tau^{3}\lVert e^{\tau\phi}u \rVert^{2}_{L^{2}(\partial\Omega)} + \tau\lVert e^{\tau\phi}\nabla u \rVert^{2}_{L^{2}(\partial\Omega)} \right) \\
         & \leq \frac{C}{\eta^{2}}\left( \lVert e^{\tau\phi}\nabla u \rVert^{2}_{L^{2}(\Omega_{\eta, 1})} + \frac{1}{\eta^{2}}\lVert e^{\tau\phi}u \rVert^{2}_{L^{2}(\Omega_{\eta, 1})} \right) + \tau^{3}\lVert e^{\tau\phi}u \rVert^{2}_{L^{2}(\partial\Omega)} + \tau\lVert e^{\tau\phi}\nabla u \rVert^{2}_{L^{2}(\partial\Omega)} .
    \end{align*}

    Choose $\tau\geq \frac{4C}{\eta^{2}}$ large enough to absorb the $H^{1}$ norm terms of $u$ over $\Omega_{\eta, 2}$ and shrink the integral domain for the left-hand side term above from $\Omega_{\eta, 2}$ to $\Omega_{\eta, 3}$ we obtain:
    \begin{align*}
         & \tau\lVert e^{\tau\phi} u\rVert^{2}_{L^{2}(\Omega_{\eta, 3})} + \lVert e^{\tau\phi} \nabla u\rVert^{2}_{L^{2}(\Omega_{\eta, 3})} \\
         & \leq \lVert e^{\tau\phi} \nabla u \rVert^{2}_{L^{2}(\Omega_{\eta, 1}\backslash
        \overline{\Omega_{\eta, 2}})} + \frac{1}{\eta^{2}}\lVert e^{\tau\phi}u \rVert^{2}_{L^{2}(\Omega_{\eta, 1}\backslash\overline{\Omega_{\eta, 2}})} + 2\tau^{2}\lVert e^{\tau\phi}u \rVert^{2}_{L^{2}(\partial\Omega)} + 2\lVert e^{\tau\phi}\nabla u \rVert^{2}_{L^{2}(\partial\Omega)}.
    \end{align*}

    Consider the maxima and minima of $\phi$ in domain $\Omega_{\eta, 3}, \Omega_{\eta, 1}\backslash\overline{\Omega_{\eta, 2}}, \partial\Omega,$ and divide both sides of the above inequality by $e^{18 \eta^{2} \tau}$ we have
    \begin{align*}
         & \lVert u\rVert^{2}_{L^{2}(\Omega_{\eta, 3})} + \lVert \nabla u\rVert^{2}_{L^{2}(\Omega_{\eta, 3})}                                                                                                                                                                                                                                                                     \\
         & \leq e^{ -10\eta^{2} \tau}\eta^{-2}\left( \lVert \nabla u \rVert^{2}_{L^{2}(\Omega_{\eta, 1}\backslash\overline{\Omega_{\eta, 2}})} + \lVert u \rVert^{2}_{L^{2}(\Omega_{\eta, 1}\backslash\overline{\Omega_{\eta, 2}})} \right) + 2\tau^{2}e^{2\tau(R_{0}^{2} -9\eta^{2}) }\left(\lVert u \rVert^{2}_{L^{2}(\partial\Omega)} + \lVert \nabla u \rVert^{2}_{L^{2}(\partial\Omega)} \right).
    \end{align*}

    We have
    \begin{align}\label{estimate u H1 M F}
        \lVert u \rVert^{2}_{H^{1}(\Omega_{\eta, 3}))} 
        \leq e^{ -10\eta^{2} \tau} \eta^{-2}\lVert u \rVert^{2}_{H^{1}(\Omega_{\eta,1})} + 2\tau^{2}e^{2\tau(R_{0}^{2} -9\eta^{2}) }\left( \lVert u \rVert_{H^{1}(\partial\Omega)} + \lVert \nabla u \rVert_{H^{1}(\partial\Omega)} \right).
    \end{align}
    On the other hand we deduce from Proposition \ref{proposition 1} the following estimate
    \begin{align} \label{zero}
        \lVert u \rVert_{H^{2}(\Omega_{\eta,1})}\leq M_{u},
    \end{align}
    where $M_{u}>0$ is a constant that only depends on $q, \Omega, k,$ and $\mathfrak p$.
    Hence, we deduce from estimate \eqref{estimate u H1 M F} and $\lVert u \rVert_{H^{1}(\partial\Omega)} + \lVert \nabla u \rVert_{H^{1}(\partial\Omega)} \leq 2\varepsilon$ as well as $\tau \geq 2$ that
    \begin{align}\label{estimate u H1 Mu epsilon}
        \lVert u \rVert^{2}_{H^{1}(\Omega_{\eta, 3}))} \leq e^{ -10\eta^{2} \tau} \left( \frac{M_{u}}{\eta}\right)^{2} + e^{2\tau(2+R_{0}^{2} -9\eta^{2}) }\varepsilon^{2}.
    \end{align}

    Finally, by assuming that $\varepsilon < \varepsilon_{0},$ with
    \begin{align}\label{varepsilon0}
        \varepsilon_{0} = e^{-(2+R_{0}^{2}-4\eta^{2})\tau_{1}}\frac{M_{u}}{\eta},
    \end{align} 
    where $\tau_{1}=\max\{\tau_{0},\frac{4C}{\eta^{2}} \}$, and $\tau_{0}$ is defined in Proposition \ref{Carleman estimate}. Let $\tau = \frac{1}{ 2+R_{0}^{2}-4\eta^{2}} \ln(\frac{M_{u}}{\eta\varepsilon}) $ and we obtain

    \begin{align*}
    \lVert u \rVert_{H^{1}(\Omega_{\eta, 3}) }\leq \sqrt{2} \left(\frac{M_{u}}{\eta}\right)^{1-\tilde \theta} \varepsilon^{\tilde \theta},
    \end{align*}
    where
    \begin{align*}
    \tilde \theta = \frac{5\eta^{2}}{2+R_{0}^{2}-4\eta^{2}}.
    \end{align*}

    By assuming that $\varepsilon\geq \varepsilon_{0}$, which implies $M_{u}\leq  e^{(2+R_{0}^{2}-4\eta^{2})\tau_{1}}\varepsilon$, choose $\tau=\tau_{1}$ in estimate \eqref{estimate u H1 Mu epsilon} and a direct consequence is

    \begin{align*}
    \lVert u \rVert_{H^{1}(\Omega_{\eta, 3}) }
    \leq \sqrt{2} e^{(2+R_{0}^{2}-9\eta^{2})\tau_{1}}\varepsilon,
    \end{align*}

    which finishes the proof.
\end{proof}

We are now ready to prove Theorem \ref{thm: stability for multiple a z}.
\begin{proof}
We treat differently the point sources with indexes within the sets $\mathfrak N, $ $\mathfrak N_{1}, $ and
$\mathfrak N_{2}$.

Assume that $i\in \mathfrak N$. We deduce from Lemma \ref{lemma Geometry}, and the fact that
$\lvert z_{i,j}-z_{i,k} \rvert \geq 8\eta$ for $i=1, 2, $ and $1\leq j\neq k \leq N_{i}$ that $B_{\eta,3}(z_{1,i})$
contains exactly two point sources $z_{1,i},z_{2,\pi(i)}$.

Lemma \ref{lemma:Carleman Estimate} and trace Theorem lead to

\begin{align}\label{u Dirichlet 1/2}
\lVert u \rVert_{H^{\frac{1}{2}}(\partial B_{3\eta}(z_{1,i}))}
\leq C_{\eta} \left( \frac{M_{u}}{\eta} \right)^{1-{\tilde \theta}} \varepsilon^{\tilde \theta},
\end{align}

if $ \varepsilon < \varepsilon_{0}, $
with $C_{\eta}>0$ is a constant that only depends on $\eta>0$.

In the domain $\Omega_{\eta,3}$, considering the equation $\Delta u + k^{2}(1+q)u=0$ and a test function $v\in H^{\frac{1}{2}}(\partial\Omega_{\eta,3})$ which can be extended as a function $v\in H^{1}(\Omega_{\eta,3})$ with extension bounded, we have 
$$
\int_{\partial\Omega_{\eta,3}} v\frac{\partial u}{\partial \nu} dS = \int_{\Omega_{\eta,3}} \nabla u \nabla v + k^{2}(1+q)u v dx.
$$

By duality of $H^{-\frac{1}{2}},H^{\frac{1}{2}},$ we have
\begin{align}\label{u Neumann -1/2}
    \lVert \frac{\partial u}{\partial\nu} \rVert_{H^{-\frac{1}{2}}(\partial B_{3\eta}(z_{1,i}))}
    \leq C_{q} \lVert u \rVert_{H^{1}(B_{3\eta}(z_{1,i}))},
\end{align}
where $C_{q}$ is a constant that depends on $k,q,\Omega,\eta$.

Combing \eqref{u Dirichlet 1/2} and \eqref{u Neumann -1/2}, we have 
\begin{align*}
\lVert  u \rVert_{L^{2}(\partial B_{3\eta}(z_{1, i}))} + \lVert \frac{\partial u}{\partial \nu} \rVert_{L^{2}(\partial B_{3\eta}(z_{1, i}))} 
\leq C_{\eta,q} \left( \frac{M_{u}}{\eta}\right)^{1-\tilde \theta} \varepsilon^{\tilde \theta},
\end{align*}
if $ \varepsilon < \varepsilon_{0},$ where $C_{\eta,q}$ is a constant that depends on $k,q,\Omega,\eta$.

Applying the results of Theorem \ref{thm: stability a z} with $S_{1} = a_{1, i} \delta_{z_{1, i}}$, $S_{2} = a_{2,\pi(i)} \delta_{z_{2, \pi(i)}}$, and $\Omega = B_{3\eta}(z_{1, i})$, we obtain the desired estimates with $C_{\mathfrak{a}}=C_{\eta,q}C_{a}M_{u}^{1-\theta},C_{\mathfrak{z}}=C_{\eta,q}C_{z}M_{u}^{1-\theta}$ and 

\begin{align} \label{theta}
\theta =  \frac{5\eta^{2}}{2+R_{0}^{2}-4\eta^{2}}.
\end{align}

if $ \varepsilon \geq \varepsilon_{0}$, analogously, choosing $\tau=\tau_{1}$ we obtain
\begin{align} \label{six}
\lVert  u \rVert_{L^{2}(\partial B_{3\eta}(z_{1, i}))} + \lVert \frac{\partial u}{\partial \nu} \rVert_{L^{2}(\partial B_{3\eta}(z_{1, i}))}  \leq C_{\eta,q}e^{(2+R_{0}^{2}-9\eta^{2})\tau_{1}} \varepsilon,
\end{align}

and the desired estimates with $\theta = 1$.

The cases where $ i\in \mathfrak N_{1}$ and $ j\in \mathfrak N_{2}$ can be treated in the same manner.
\end{proof}

\begin{remark}
    If $ \varepsilon < \varepsilon_{0}$, it follows from the dependence of the constants with respect to $\eta$ in the above proof that $M_{u}$ grows polynomially with order $-1$ by combining the results of \cite{gruter1982green} while $\varepsilon_{0}$ given explicitly in \eqref{varepsilon0} decays exponentially as $\eta$ approaches zero. The constant factor $C_{\eta,q}$ of $ C_{\mathfrak{a}} $ and $C_{\mathfrak{z}}$ blows up polynomially with negative order following the results of \cite{fernandezbonder2002asymptotic} while $C_{a},C_{z}$ decays polynomially with order $1$ when $\eta$ approaches zero. This means that when the point sources are closer, it becomes more difficult to recover them.
\end{remark}

\section{Conclusion}

In this paper, we have derived a H\"{o}lder type stability estimate for the inverse problem of recovering a finite number of point sources. The prior separation information of point sources turns out to be crucial in the proof of this result. A future project is to explore whether the geometric assumption could be removed. We also plan to study the increasing stability of the inversion as the frequency $k$ grows \cite{bao2010error, bao2010multifrequency, hrycak2004increased}.
\appendix

\section{Complex geometrical optics solution}\label{section: CGO solution}
We restrict the Helmholtz equation to the cube $D:=(-R_{0}, R_{0})^{3}\subset \mathbb{R}^{3}$ to apply the periodic Faddeev-type operator which is an important tool in the reconstruction of CGO solutions \cite{hahner1996periodic}. The shifted grid $\Gamma$, and the corresponding orthonormal basis $\{e_{\alpha}\}$ in $L^{2}(D)$ for the Fourier series are introduced as
$$
    \Gamma:=
    \left\{
    \alpha=(\alpha_{1},\alpha_{2},\alpha_{3})^{T}\in \mathbb{R}^{3}:\frac{R_{0}}{\pi}\alpha_{2}-\frac{1}{2}\in Z,\frac{R_{0}}{\pi}\alpha_{j}\in \mathbb{Z},j=1,3
    \right\}
$$
with
$e_{\alpha}=(2R_{0})^{-3/2}e^{i\alpha \cdot x},x\in D,\alpha \in \Gamma.$

For $f\in L^{2}(D),$ we have $f=\sum_{\alpha \in \Gamma}\hat{f} e_{\alpha}$ where $\hat{f}=\int_D f \overline{e_{\alpha}} dx$. We present the following result for the periodic solution $u=G_{\xi}f$ of $(\Delta	+ 2i\xi \cdot \nabla) u = f$.
\begin{lemma}\cite{hahner1996periodic}\label{lemma:faddeev type op.}
    Let $s\in \mathbb{R}, t>0$ be a real number, and $\xi:=(s,it,0)^{T}\in \mathbb{C}^{3}.$ Then, the operator
    $$
        G_{\xi}:L^{2}(D)\longrightarrow H^{2}(D) \quad
        G_{\xi} f:=-\sum_{\alpha \in \Gamma}\frac{\hat{f}(\alpha)}{(\alpha\cdot\alpha + 2\xi\cdot\alpha)}e_{\alpha},
    $$
    is well defined, which satisfies
    \begin{align*}
         & \lVert G_{\xi}f \rVert_{L^{2}(D)} \leq \frac{R_{0}}{\pi t} \lVert f\rVert_{L^{2}(D)}, \\
         & \sum_{i=1}^{3}\lVert \partial_{x_{i}} G_{\xi}f \rVert_{L^{2}(D)}
        \leq \frac{R_{0}}{\pi }\frac{\lvert s\rvert+\sqrt{\lvert s\rvert^{2}+\pi t /R_{0}}}{t} \lVert f\rVert_{L^{2}(D)},
    \end{align*}
    for all $f\in L^{2}(D)$. In addition $(\Delta + 2i\xi \cdot \nabla)G_{\xi}f = f $ in the weak sense for all $f\in L^{2}(D)$.
\end{lemma}


Recursively using Lemma \ref{lemma:faddeev type op.}, we can reproduce the existence of a periodic CGO solution and a decaying property of its associated potential in $C^{0}(\overline D)$ as $\lvert \Im\xi\rvert$ becomes large.
\begin{lemma}\label{lemma:cgo decay n2}
    Let $q\in C^0(\overline D)$, $\xi \in \mathbb{C}^{3},\xi\cdot\xi=k^{2}$, and
    $\frac{2R_{0} k^{2}}{\pi \lvert \Im \xi\rvert} \lVert q \rVert_{C^{0}(\overline D)}\leq 1$, then the Helmholtz equation \eqref{equation: forward problem} has a periodic solution $u$
    in $D$, with the following form
    $$u=e^{ix\cdot\xi}(1+\phi_{\xi}(x)),$$
    and
    $$
        \lVert \phi_{\xi}\rVert_{L^{2}(D)}
        \leq \frac{2R_{0}k^{2}}{\pi\lvert\Im\xi\rvert}\lVert q\rVert_{L^{2}(D)}.
    $$

    Furthermore,
    if $q\in C^{2}(\overline D)$ and $\max\{\lVert q\rVert_{C^{2}(\overline D)} \frac{2R_{0} k^{2}}{\pi},\frac{\pi}{R_{0}},k^{2}\} \leq \lvert \Im \xi \rvert$, then
    \begin{align*}
         & \lVert \phi_{\xi}\rVert_{H^{2}(D)}
        \leq \frac{243 R_{0}k^{2}}{\pi\lvert\Im\xi\rvert}\lVert q \rVert_{H^{2}(D)}, \\
         & \lVert \phi_{\xi}\rVert_{H^{3}(D)}
        \leq M_{G},
    \end{align*}
    where $M_{G}=\left(\left(\frac{R_{0}k^{2}}{\pi}(4\lVert q\rVert_{C^{0}(\overline D)}+24)+\frac{27}{2}\right)
    \frac{18R_{0}k^{2}}{\pi |\Im\xi|} + \frac{4 R_{0}k^{2}}{\pi} \right)\lVert q \rVert_{H^{2}(D)}.$
\end{lemma}

\begin{proof}
    The operator $\Delta$ is invariant under an unitary coordinate transform, we can assume that $$\xi=(s,it,0)^{T} \in \mathbb{C}^{3}, \xi\cdot\xi=k^{2}.$$

    A direct calculation shows that $\phi_{\xi}$ satisfies
    \begin{align}\label{equation:faddeev for phi}
        (\Delta	+ 2i\xi \cdot \nabla)\phi_{\xi} + k^{2}q\phi_{\xi}=-k^{2}q.
    \end{align}

    From Lemma \ref{lemma:faddeev type op.} we deduce that $G_{\xi}=(\bigtriangleup+2i\xi\cdot\nabla)^{-1}$ is a bounded operator on $L^{2}(D)$ with norm bounded by $\frac{R_{0}}{\pi t}$. Applying $G_{\xi}$ on both sides of \eqref{equation:faddeev for phi}, we have
    $$
        \phi_{\xi} + G_{\xi}\left(k^{2}q\phi_{\xi}\right)=G_{\xi}(-k^{2}q).
    $$

    Since $q\in C^0(\overline D)$, $G_{\xi}(k^{2}q\cdot )$ is a bounded operator on $L^{2}(D)$
    satisfying

    $$ 
    \|G_{\xi}(k^{2}q\cdot )\|_{\mathcal L(L^{2}(D))} \leq k^{2}\lVert q \rVert_{C^{0}(\overline D)}\lVert G_{\xi}\rVert_{\mathcal L(L^{2}(D))}.
    $$

    Lemma \ref{lemma:faddeev type op.} and inequality $\frac{2R_{0} k^{2}}{\pi \lvert \Im \xi\rvert}\lVert q \rVert_{C^{0}(\overline D)}\leq 1$, lead to

    $$ \|G_{\xi}(k^{2}q\cdot )\|_{\mathcal L(L^{2}(D))} \leq\frac{1}{2}.$$

    Hence $I+G_{\xi}(k^{2}q(x))$ is invertible with $\lVert \left(I+G_{\xi}(k^{2}q)\right)^{-1}\rVert_{\mathcal L(L^{2}(D))} \leq 2$. Consequently

    $$
        \lVert \phi_{\xi} \rVert_{L^{2}(D)}
        =\lVert (I+G_{\xi}(k^{2}q))^{-1}G_{\xi}(-k^{2}q)\rVert_{L^{2}(D)}
        \leq 2\lVert G_{\xi}(-k^{2}q)\rVert_{L^{2}(D)}
        \leq \frac{2R_{0}k^{2}}{\pi t}\lVert q \rVert_{L^{2}(D)}.
    $$
    Next, we estimate the $H^{s}$ norm of $\phi_{\xi}$. Taking the derivative with respect to $x_{i}$ on both sides of \eqref{equation:faddeev for phi}, we obtain
    $$
        (\Delta + 2i\xi \cdot \nabla)\partial_{x_{i}}\phi_{\xi} + k^{2}q\partial_{x_{i}}\phi_{\xi}
        =-k^{2}\partial_{x_{i}}q-k^{2}\partial_{x_{i}}q\phi_{\xi}.
    $$

    Applying again $G_{\xi}$ on both sides of the equation above, with the assumption that $\lVert q\rVert_{C^{2}(\overline D)} \frac{2R_{0} k^{2}}{\pi \lvert \Im\xi \rvert} \leq 1$, and the decay estimate of $\phi_{\xi}$, we determine
    the estimate of $\partial_{x_{i}}\phi_{\xi}$ as
    \begin{align*}
        \lVert \partial_{x_{i}}\phi_{\xi} \rVert_{L^{2}(D)}
        \leq \frac{2R_{0}k^{2}}{\pi t}\lVert \partial_{x_{i}}q+ \partial_{x_{i}}q\phi_{\xi}\rVert_{L^{2}(D)}
        \leq\frac{4R_{0}k^{2}}{\pi t} \lVert q\rVert_{H^{1}(D)}.
    \end{align*}

    Similarly, considering the second-order derivatives of $\phi$, we have
    $$
        (\Delta	+ 2i\xi \cdot \nabla)\partial_{x_{i}}\partial_{x_{j}}\phi_{\xi}
        + k^{2}q(x)\partial_{x_{i}}\partial_{x_{j}}\phi_{\xi}
        = -k^{2}\partial_{x_{i}}\partial_{x_{j}}q(x)(1+\phi_{\xi})
        -k^{2}\partial_{x_{i}}q\partial_{x_{j}}\phi_{\xi}
        -k^{2}\partial_{x_{j}}q\partial_{x_{i}}\phi_{\xi}, \; i,j =1, 2, 3.
    $$

    Applying $(1+G_{\xi}k^{2}q)^{-1}G_{\xi}$ on both sides of the above equation and estimate the four terms on the right hand side side separately with the assumption that $\lVert q\rVert_{C^{2}(\overline D)} \frac{2R_{0} k^{2}}{\pi \lvert \Im \xi\rvert}\leq 1$, we obtain
    $$
        \lVert\partial_{x_{i}}\partial_{x_{j}}\phi_{\xi}\rVert_{L^{2}(D)}
        \leq
        \frac{2R_{0}k^{2}}{\pi t}
        \lVert \partial_{x_{i}}\partial_{x_{j}}q \rVert_{L^{2}(D)}
        + \lVert \phi_{\xi}\rVert_{L^{2}(D)}
        + \lVert \partial_{x_{i}}\phi_{\xi}\rVert_{L^{2}(D)}
        + \lVert \partial_{x_{j}}\phi_{\xi}\rVert_{L^{2}(D)}.
    $$
    Using the previous results for $\phi_{\xi}$ and $\partial_{i}\phi_{\xi},\partial_{j}\phi_{\xi}$, we have
    \begin{align*}
        \lVert\partial_{x_{i}}\partial_{x_{j}}\phi_{\xi}\rVert_{L^{2}(D)}
         & \leq
        \frac{2R_{0}k^{2}}{\pi t}
        \left(
        \lVert \partial_{x_{i}}\partial_{x_{j}}q \rVert_{L^{2}(D)}
        + C_{2}^{0}\lVert q \rVert_{H^{0}(D)}
        + 2C_{2}^{1}\lVert q \rVert_{H^{1}(D)}
        \right) \leq
        \frac{18R_{0}k^{2}}{\pi t}\lVert q \rVert_{H^{2}(D)}.
    \end{align*}

    For the third-order derivative of $\phi_{\xi}$, concerning the second estimate of Lemma \ref{lemma:faddeev type op.} and the following equation
    $$
        (\Delta	+ 2i\xi \cdot \nabla)\partial_{x_{i}}\phi_{\xi}
        = - k^{2}q\partial_{x_{i}}\phi_{\xi} -k^{2}\partial_{x_{i}}q-k^{2}\partial_{x_{i}}q\phi_{\xi},
    $$
    we have
    \begin{align*}
        \sum_{l=1}^{3}\lVert \partial^{3}_{x_{l}x_{i}x_{j}}\phi_{\xi} \rVert_{L^{2}(D)}
         & \leq
        \frac{R_{0}k^{2}}{\pi}\frac{\lvert s \rvert + \sqrt{ \lvert s \rvert^{2} + \frac{\pi t}{R_{0}}}}{t}
        \lVert
        q\partial^{2}_{x_{i}x_{j}}\phi_{\xi}
        \rVert_{L^{2}(D)} \\
         & +
        \frac{R_{0}k^{2}}{\pi}\frac{\lvert s \rvert + \sqrt{ \lvert s \rvert^{2} + \frac{\pi t}{R_{0}}}}{t}
        \lVert
        \partial^{2}_{x_{i}x_{j}} q (1+\phi_{\xi}) + \partial_{x_{i}}q\partial_{x_{j}}\phi_{\xi} + \partial_{x_{j}}q \partial_{x_{i}}\phi_{\xi}
        \rVert            \\
         & \leq
        \frac{R_{0}k^{2}}{\pi}\frac{\lvert s \rvert + \sqrt{ \lvert s \rvert^{2} + \frac{\pi t}{R_{0}}}}{t}
        \left((\lVert q\rVert_{C^{0}(\overline D)}+6)\frac{18R_{0}k^{2}}{\pi t}+1\right)\lVert q \rVert_{H^{2}(D)}.
    \end{align*}

    Notice that $\frac{\pi}{R_{0}},k^{2} \leq \lvert \Im \xi \rvert = t, s^{2}=t^{2}+k^{2}$, a direct calculation then gives
    $$
        \sum_{l=1}^{3}\lVert \partial_{x_{l}}\partial_{x_{i}}\partial_{x_{j}}\phi_{\xi} \rVert_{L^{2}(D)}
        \leq
        \frac{R_{0}k^{2}}{\pi}(\sqrt{2}+\sqrt{3})
        \left((\lVert q\rVert_{C^{0}(D)}+6)\frac{18R_{0}k^{2}}{\pi t}+1\right) \lVert q \rVert_{H^{2}(D)}.
    $$

    Collecting all terms and via a rough estimate of the constant we find
    \begin{align*}
         & \lVert \phi_{\xi}\rVert_{H^{1}(D)}
        \leq \frac{18 R_{0}k^{2}}{\pi\lvert\Im\xi\rvert}\lVert q \rVert_{H^{1}(D)}, \\
         & \lVert \phi_{\xi}\rVert_{H^{2}(D)}
        \leq \frac{243 R_{0}k^{2}}{\pi\lvert\Im\xi\rvert}\lVert q \rVert_{H^{2}(D)}, \\
         & \lVert \phi_{\xi}\rVert_{H^{3}(D)}
        \leq
        \left(\left(\frac{R_{0}k^{2}}{\pi}(4\lVert q\rVert_{C^{0}(\overline D)}+24)+\frac{27}{2}\right)
        \frac{18R_{0}k^{2}}{\pi |\Im\xi|} + \frac{4 R_{0}k^{2}}{\pi} \right)\lVert q \rVert_{H^{2}(D)}.
    \end{align*}
    Finally, by noticing $\max\{\lVert q\rVert_{C^{2}(\overline D)} \frac{2R_{0} k^{2}}{\pi},\frac{\pi}{R_{0}}, k^{2}\} \leq \lvert \Im \xi \rvert$, we obtain the desired result.
\end{proof}

\begin{remark}\label{remark:1/12 gradient decay and continous dependence on xi}
    In Lemma \ref{lemma:cgo decay n2}, if in addition, we assume $\lvert \Im \xi\rvert \geq \frac{2916 R_{0}k^{2}}{\pi}\lVert q(x) \rVert_{H^{2}(D)}C_{L}$ where $C_{L}>0$ is the constant of the continuous embeding $H^{2}(D)\hookrightarrow C^{0}(\overline{D})$ depending only on $D$, then
    \begin{align*}
        \lVert \phi_{\xi}\rVert_{C^{0}(\overline D)}\leq C_{L} \lVert \phi_{\xi}\rVert_{H^{2}(D)}\leq \frac{1}{12},\\
        \lVert \partial_{x_{i}}\phi_{\xi}\rVert_{C^{0}(\overline D)} \leq \lVert \partial_{x_{i}}\phi_{\xi}\rVert_{H^{2}(D)} \leq M_{G}C_{L}.
    \end{align*}
    Consequently
    \begin{align*}
        \lvert v_{\xi}(x)\rvert\leq 2 e^{tR_{0}}, \lvert \nabla v_{\xi}(x) \rvert\leq e^{tR_{0}}(4t+2k +M_{G}C_{L}), \quad \forall x\in D.
    \end{align*}
\end{remark}

The following part is about the continuity of the CGO potential $\phi_{\xi}$ on the parameter $\xi$ which will be used in Theorem \ref{thm: stability a z}.

\begin{proposition} \label{Holder continuity potential}
    For a fixed point $x_{0} \in \mathbb{C}^{3},$ let $\phi_{\xi}$ be the potential of the CGO solution constructed as in Lemma \ref{lemma:cgo decay n2}, then
    $$
        \phi(x_{0}): \mathbb{C}^{3}\longrightarrow \mathbb{C}, \; \;
        \xi \rightarrow \phi_{\xi}(x_{0}),
    $$
    is a well-defined continuous function on the complex set $\lvert \Im \xi\rvert \geq \max\{\lVert q\rVert_{C^{2}(\overline D)} \frac{2R_{0} k^{2}}{\pi},\frac{\pi}{R_{0}},k^{2}\}$.

\end{proposition}
\begin{proof}
The result can be verified based on the following observation. Considering the equation of $\phi_{\xi_{1}}-\phi_{\xi_{2}}$ as below,
$$
    (\bigtriangleup + 2i\xi_{1} \cdot \nabla)(\phi_{\xi_{1}}-\phi_{\xi_{2}}) + k^{2}(1+q(x))(\phi_{\xi_{1}}-\phi_{\xi_{2}})=(\xi_{1}-\xi_{2})\cdot\nabla\phi_{\xi_{2}},
$$
it follows the same idea in the proof of Lemma \ref{lemma:cgo decay n2} and Remark \ref{remark:1/12 gradient decay and continous dependence on xi} that
\begin{align*}
    \lvert \phi(x_{0})(\xi_{1}) - \phi(x_{0})(\xi_{2}) \rvert
     & \leq C_{L} \lVert \phi_{\xi_{1}} - \phi_{\xi_{2}} \rVert_{H^{2}(D)}
    \leq 52C_{L} \left(1+k^{2}\lVert q \rVert_{C^{2}(\overline{D})} \right)^{2}\lVert \phi_{\xi_{2}}\rVert_{H^{3}(D)}\lvert \xi_{1}-\xi_{2}\rvert,
\end{align*}
which shows the continuity for $\phi(x_{0})$.

\end{proof}

\section{Carleman estimates}

Let $\phi \in C^{2}(\overline{\omega})$ without critical points in $\omega$ and $\varphi =e^{\lambda \phi}$. Recall the Carleman estimate for the Helmholtz equation \cite{hrycak2004increased, choulli2016applications, klibanov1992inverse}.

\begin{proposition}[\cite{hrycak2004increased}] \label{Carleman estimate} Let $\omega \subset B_{\frac{R_{0}}{2}}$ be an open connected domain satisfying $\omega\cap B_{\epsilon}(z_{0})=0$ for some $z_{0}\in \mathbb{R}^{3}$ and a constant $0 < \epsilon \in \mathbb{R}$ small enough. Let $\phi=\lvert x-z_{0} \rvert^{2}$, there exist three positive constants $C$, and $\tau_{0}$, which depend only on
    $\epsilon$, $k$, $R_{0}$, and $\|q\|_{C(\overline \omega)}$ so that
    \begin{align*}
         & \tau^{2}\lVert e^{\tau\phi} u \rVert^{2}_{L^{2}(\omega)} + \tau\lVert e^{\tau\phi} \nabla u\rVert^{2}_{L^{2}(\omega)}                                                                                          \\
         & \leq C \left( \lVert e^{\tau\phi}A u \rVert^{2}_{L^{2}(\omega)} + \tau^{3}\lVert e^{\tau\phi}u \rVert^{2}_{L^{2}(\partial\omega)} + \tau\lVert e^{\tau\phi}\nabla u \rVert^{2}_{L^{2}(\partial\omega)} \right).
    \end{align*}

    for all $u\in H^{2}(\omega),$ and $\tau\ge\tau_{0}=\max\{2,\frac{k^{2}(1+\lvert q(x)\rvert_{C(\overline\omega)})}{\epsilon},\frac{k^{2}\lvert \nabla q(x)\rvert_{C(\overline\omega)}}{\epsilon},\frac{1}{\epsilon^{2}}\}$, where $A u=\Delta u + k^{2}(1+q(x))u $.
\end{proposition}
\begin{proof}
    The main idea in this proof is from Lemma 2.2 of \cite{hrycak2004increased}. Here, we give the details to derive an explicit expression for the constants $C,\tau_{0}$.

    From the weight function $\phi=\lvert x-z_{0} \rvert^{2},$ a direct calculation gives $\nabla \phi= 2(x-z_{0}),\nabla\phi\cdot\nabla\phi=4\lvert x-z_{0} \rvert^{2},\Delta\phi=6$.

    Let $v=e^{\tau \phi} u,$ and consider 
    \begin{align*}
        &e^{\tau \phi} A u=e^{\tau \phi}A(e^{-\tau\phi}v)=-\tau v\Delta\phi + \tau^{2}v\nabla\phi\cdot\nabla\phi-2\tau\nabla\phi\cdot\nabla v + k^{2}(1+q(x))v+\Delta v\\
        &=-6\tau v - 4\tau (x-z_{0})\cdot\nabla v +\left(4\tau^{2}(x-z_{0})^{2}+k^{2}(1+q(x))\right)v+\Delta v.
    \end{align*} 

    Then
    \begin{align*}
        &\int_{\omega} (e^{\tau \phi} A u)^{2} d x \\
        &\geq  \int_{\omega} \left(-6\tau v - 4\tau (x-z_{0})\cdot\nabla v +\left(4\tau^{2}(x-z_{0})^{2}+k^{2}(1+q(x))\right)v+\Delta v \right)^{2}d x\\
        &-\int_{\omega} \left( 6\tau v + 4\tau (x-z_{0})\cdot\nabla v +\left(4\tau^{2}(x-z_{0})^{2}+k^{2}(1+q(x))\right)v+\Delta v \right)^{2}d x\\
        &=\int_{\omega}4\left(6\tau v + 4\tau (x-z_{0})\cdot\nabla v\right)\Bigl(\bigl(4\tau^{2}(x-z_{0})^{2}+k^{2}(1+q(x))\bigr)v+\Delta v \Bigr)d x\\
        &= 4\int_{\omega} (-6\tau v \left(4\tau^{2}(x-z_{0})^{2}+k^{2}(1+q(x))\right)v) d x + 4\int_{\omega} -6\tau v \Delta v d x\\
        & + 4\int_{\omega} -4\tau(x-z_{0})\cdot\nabla v \left(4\tau^{2}(x-z_{0})^{2}+k^{2}(1+q(x))\right)v d x + 4\int_{\omega} -4\tau(x-z_{0})\cdot\nabla v \Delta v d x\\
        &= 4(I_{1}+I_{2}+I_{3}+I_{4}).
    \end{align*}
    Here 
    $$I_{1}=-6\tau \int_{\omega} \Bigl(4\tau^{2}(x-z_{0})^{2}+k^{2}\bigl(1+q(x)\bigr)\Bigr)v^{2}  d x,$$
    \begin{align*}
        I_{2}=-6\tau \int_{\omega}  v \Delta v d x = -6\tau \int_{\partial\omega} v\frac{\partial v}{\partial\nu} d S + 6\tau \int_{\omega}\lvert \nabla v\rvert^{2} d x,
    \end{align*}
    \begin{align*}
        &I_{3} = - 4\tau\int_{\omega}  (x-z_{0})\cdot \nabla v \Bigl(\bigl(4\tau^{2}(x-z_{0})^{2}+k^{2}(1+q(x))\bigr)v \Bigr) d x\\
        &= - 2\tau\int_{\omega}  (x-z_{0})\cdot \nabla v^{2} \bigl(4\tau^{2}(x-z_{0})^{2}+k^{2}(1+q(x))\bigr)  d x\\
        &=- 2\tau\int_{\partial\omega} (x-z_{0})\cdot \nu \bigl(4\tau^{2}(x-z_{0})^{2}+k^{2}(1+q(x))\bigr) v^{2} d S \\
        &  + 2\tau\int_{\omega}\mbox{div}\Bigl( (x-z_{0})\bigl(4\tau^{2}(x-z_{0})^{2}+k^{2}(1+q(x))\bigr) \Bigr) v^{2} d x \\
        &= - 2\tau\int_{\partial\omega} (x-z_{0})\cdot \nu \bigl(4\tau^{2}(x-z_{0})^{2}+k^{2}(1+q(x))\bigr) v^{2} d S\\
        &  + 2 \tau \int_{\omega} \left(20\tau^{2}\lvert x-z_{0}\rvert^{2} + 3 k^{2}(1+q(x))\right) v^{2} +  k^{2}\nabla q(x)\cdot (x-z_{0}) v^{2}  d x,
    \end{align*}
    \begin{align*}
        &I_{4} = - 4 \tau \int_{\omega} (x-z_{0})\cdot \nabla v \Delta v d x\\
        &=  -4 \tau \int_{\partial\omega} (x-z_{0})\cdot\nabla v \frac{\partial v}{\partial \nu} d S + 4\tau \int_{\omega} \nabla \left( (x-z_{0})\cdot \nabla v\right)\cdot \nabla v\\
        &=  -4 \tau \int_{\partial\omega} (x-z_{0})\cdot\nabla v \frac{\partial v}{\partial \nu} d S + 4\tau\int_{\omega} \lvert \nabla v \rvert^{2} + (x-z_{0})\frac{1}{2}\nabla\lvert \nabla v\rvert^{2} d x\\
        &=  -4 \tau \int_{\partial\omega} (x-z_{0})\cdot\nabla v \frac{\partial v}{\partial \nu} d S + 4\tau\int_{\omega} \lvert \nabla v \rvert^{2} d x + 2\tau \int_{\partial\omega}\nu\cdot (x-z_{0}) \lvert \nabla v\rvert^{2} d S - 4\tau \int_{\omega} \frac{3}{2} \lvert \nabla v\rvert^{2} d x.
    \end{align*}
    Thus,
    \begin{align*}
        &I_{1}+I_{2}+I_{3}+I_{4}\\
        &\geq 4\tau \int_{\omega}  \lvert \nabla v \rvert^{2} + \left(4\tau^{2}\lvert x-z_{0} \rvert^{2}  + \frac{1}{2} k^{2}\nabla q(x) (x-z_{0}) \right) v^{2}  d x\\
        &-6\tau \int_{\partial\omega}   v\frac{\partial v}{\partial \nu} d S - 2\tau\int_{\partial\omega}\left( 4\tau^{2}\lvert x-z_{0} \rvert^{2} + k^{2}(1+q(x)) (x-z_{0})\cdot\nu \right) v^{2} d S \\
        &- 4\tau \int_{\partial\omega} (x-z_{0})\cdot\nabla v \frac{\partial v}{\partial \nu} d S + 2\tau\int_{\partial \omega} \nu\cdot(x-z_{0}) \lvert \nabla v \rvert^{2} d S.
    \end{align*}
    Notice that $\tau\ge\tau_{0}=\max\{2,\frac{k^{2}(1+\lvert q(x)\rvert_{C(\overline\omega)})}{\epsilon},\frac{k^{2}\lvert \nabla q(x)\rvert_{C(\overline\omega)}}{\epsilon},\frac{1}{\epsilon^{2}}\},\lvert x-z_{0}\rvert \geq \epsilon$, then in the domain $\omega$ we have
    $$
    \tau \left(k^{2}(1+q(x)) (x-z_{0})\cdot\nu \right) v^{2},\tau\left(k^{2}\nabla q(x) (x-z_{0})\right)v^{2}\leq \tau^{2}\lvert x-z_{0} \rvert^{2}v^{2}.
    $$

    Since $\omega\subset B_{\frac{R_{0}}{2}}$ we have
    \begin{align*}
    \frac{1}{4}\int_{\omega} (e^{\tau \phi} A u)^{2} d x
    &\geq 4\tau \int_{\omega}  \lvert \nabla v \rvert^{2} +  2\tau^{2}\lvert x-z_{0} \rvert^{2} v^{2}  d x\\
    &- \left(10\tau^{3}R_{0}^{2} + 3\tau \right)\int_{\partial\omega} v^{2} d S - \left(6\tau R_{0} + 3\tau \right)\int_{\partial\omega} \lvert \nabla v\rvert^{2} d S.
    \end{align*}
    Combing $ \nabla v = e^{\tau\phi}( \tau u \nabla \phi  + \nabla u) = e^{\tau\phi}( \tau (x-z_{0}) u  + \nabla u)$, $\omega\subset B_{\frac{R_{0}}{2}}$ and $\lvert x- z_{0}\rvert \geq \epsilon$, we get

    \begin{align*}
    \frac{1}{4}\int_{\omega} (e^{\tau \phi} A u)^{2} d x
    &\geq 4\tau \int_{\omega}  \frac{1}{2}\lvert e^{\tau\phi}\nabla u \rvert^{2} +  \tau^{2}\epsilon^{2}e^{2\tau\phi}u^{2}  d x\\
    &- (10\tau^{3}R_{0}^{2} +  3\tau )\int_{\partial\omega} e^{2\tau\phi} u^{2} d S - (6\tau R_{0} + 3\tau)\int_{\partial\omega} \lvert \nabla (e^{\tau\phi}u) \rvert^{2} d S.
    \end{align*}

    Notice that $\tau\epsilon^{2}\geq 1$, reformulating the above inequality we obtain the desired result, that is,

    \begin{align*}
        2\tau \int_{\omega}  \lvert e^{\tau \phi} \nabla u \rvert^{2} +  \tau(e^{\tau \phi}u)^{2}  d x  
        &\leq  \frac{1}{4} \int_{\omega} (e^{\tau \phi} A u)^{2} d x  +  (6\tau R_{0} + 3\tau)\int_{\partial\omega} \lvert e^{\tau\phi}\nabla u\rvert^{2} d S \\
        &+  \left( 10\tau^{3}R_{0}^{2} + 3\tau + 4\tau^{2} R_{0}^{2}(6\tau R_{0} + 3\tau) \right)\int_{\partial\omega} e^{2\tau\phi}u^{2} d S\\
        &\leq C \left( \int_{\omega} (e^{\tau \phi} A u)^{2} d x + \tau\int_{\partial\omega} \lvert e^{\tau\phi}\nabla u\rvert^{2} d x + \tau^{3}\int_{\partial\omega}  e^{2\tau\phi} u^{2} d x \right).
    \end{align*}
\end{proof}

\begin{remark}
    The result of Proposition \ref{Carleman estimate} holds for a genreal case $\Omega\subset \mathbb{R}^{n},$ and for further details we refer the reader to \cite{hahner1996periodic}.
\end{remark}

\bibliography{mybib}

\begin{thebibliography}{10}

\bibitem{acosta2012multifrequency}
S.~Acosta, S.~Chow, J.~Taylor, and V.~Villamizar.
\newblock On the multi-frequency inverse source problem in heterogeneous media.
\newblock {\em Inverse Problems}, 28(7):075013, 2012.

\bibitem{alessandrini2009stability}
G.~Alessandrini, L.~Rondi, E.~Rosset, and S.~Vessella.
\newblock The stability for the {{Cauchy}} problem for elliptic equations.
\newblock {\em Inverse Problems}, 25(12):123004, 2009.

\bibitem{alves2009full}
C.~J.~S. Alves, N.~F.~M. Martins, and N.~C. Roberty.
\newblock Full identification of acoustic sources with multiple frequencies and
  boundary measurements.
\newblock {\em Inverse Problems and Imaging}, 3(2):275--294, 2009.

\bibitem{ammari2002inverse}
H.~Ammari, G.~Bao, and J.~L. Fleming.
\newblock An inverse source problem for {{Maxwell}}'s equations in
  magnetoencephalography.
\newblock {\em SIAM Journal on Applied Mathematics}, 62(4):1369--1382, 2002.

\bibitem{angel1991antenna}
T.~Angel, A.~Kirsch, and R.~Kleinmann.
\newblock Antenna control and generalized characteristic modes.
\newblock {\em Proceedings of the IEEE}, 79:1559--1568, 1991.

\bibitem{bao2010multifrequency}
G.~Bao, J.~Lin, and F.~Triki.
\newblock A multi-frequency inverse source problem.
\newblock {\em Journal of Differential Equations}, 249(12):3443--3465, 2010.

\bibitem{bao2011numerical}
G.~Bao, J.~Lin, and F.~Triki.
\newblock Numerical solution of the inverse source problem for the {H}elmholtz
  equation with multiple frequency data.
\newblock In {\em Mathematical and statistical methods for imaging}, volume 548
  of {\em Contemp. Math.}, pages 45--60. Amer. Math. Soc., Providence, RI,
  2011.

\bibitem{bao2021recovering}
G.~Bao, Y.~Liu, and F.~Triki.
\newblock Recovering simultaneously a potential and a point source from
  {{Cauchy}} data.
\newblock {\em Minimax Theory and its Applications}, 6(2):227--238, 2021.

\bibitem{bao2015recursive}
G.~Bao, S.~Lu, W.~Rundell, and B.~Xu.
\newblock A {{Recursive Algorithm}} for {{MultiFrequency Acoustic Inverse
  Source Problems}}.
\newblock {\em SIAM Journal on Numerical Analysis}, 53(3):1608--1628, 2015.

\bibitem{bao2010error}
G.~Bao and F.~Triki.
\newblock Error estimates for the recursive linearization of inverse medium
  problems.
\newblock {\em Journal of Computational Mathematics}, 28(6):725--744, 2010.

\bibitem{bleistein1977nonuniqueness}
N.~Bleistein and J.~K. Cohen.
\newblock Nonuniqueness in the inverse source problem in acoustics and
  electromagnetics.
\newblock {\em Journal of Mathematical Physics}, 18(2):194, 1977.

\bibitem{choulli2016applications}
M.~Choulli.
\newblock {\em Applications of {{Elliptic Carleman Inequalities}} to {{Cauchy}}
  and {{Inverse Problems}}}.
\newblock {{SpringerBriefs}} in {{Mathematics}}. {Springer International
  Publishing}, {Cham}, 2016.

\bibitem{colton1998inverse}
D.~Colton and R.~Kress.
\newblock {\em Inverse {{Acoustic}} and {{Electromagnetic Scattering Theory}}},
  volume~93 of {\em Applied {{Mathematical Sciences}}}.
\newblock {Springer Berlin Heidelberg}, {Berlin, Heidelberg}, 1998.

\bibitem{elbadia2005inverse}
A.~El~Badia.
\newblock Inverse source problem in an anisotropic medium by boundary
  measurements.
\newblock {\em Inverse Problems}, 21(5):1487--1506, 2005.

\bibitem{elbadia2012holder}
A.~El~Badia and A.~El~Hajj.
\newblock H\"older stability estimates for some inverse pointwise source
  problems.
\newblock {\em Comptes Rendus Mathematique}, 350(23-24):1031--1035, 2012.

\bibitem{elbadia2011inverse}
A.~El~Badia and T.~Nara.
\newblock An inverse source problem for {{Helmholtz}}'s equation from the
  {{Cauchy}} data with a single wave number.
\newblock {\em Inverse Problems}, 27(10):105001, 2011.

\bibitem{fernandezbonder2002asymptotic}
J.~Fern{\'a}ndez~Bonder and J.~D. Rossi.
\newblock Asymptotic behavior of the best {{Sobolev}} trace constant in
  expanding and contracting domains.
\newblock {\em Communications on Pure and Applied Analysis}, 1(3):359--378,
  2002.

\bibitem{fokas2004unique}
A.~S. Fokas, Y.~Kurylev, and V.~Marinakis.
\newblock The unique determination of neuronal currents in the brain via
  magnetoencephalography.
\newblock {\em Inverse Problems}, 20(4):1067--1082, 2004.

\bibitem{gruter1982green}
M.~Gr{\"u}ter and K.-O. Widman.
\newblock The {{Green}} function for uniformly elliptic equations.
\newblock {\em Manuscripta Mathematica}, 37(3):303--342, 1982.

\bibitem{hahner1996periodic}
P.~H{\"a}hner.
\newblock A {{Periodic Faddeev}}-{{Type Solution Operator}}.
\newblock {\em Journal of Differential Equations}, 128(1):300--308, 1996.

\bibitem{hamalainen1993magnetoencephalography}
M.~H{\"a}m{\"a}l{\"a}inen, R.~Hari, R.~J. Ilmoniemi, J.~Knuutila, and O.~V.
  Lounasmaa.
\newblock Magnetoencephalography\textemdash theory, instrumentation, and
  applications to noninvasive studies of the working human brain.
\newblock {\em Reviews of Modern Physics}, 65(2):413--497, 1993.

\bibitem{he1998identification}
S.~He and V.~G. Romanov.
\newblock Identification of dipole sources in a bounded domain for
  {{Maxwell}}'s equations.
\newblock {\em Wave Motion}, 28(1):25--40, 1998.

\bibitem{hrycak2004increased}
T.~Hrycak and V.~Isakov.
\newblock Increased stability in the continuation of solutions to the
  {{Helmholtz}} equation.
\newblock {\em Inverse Problems}, 20(3):697--712, 2004.

\bibitem{klibanov1992inverse}
M.~V. Klibanov.
\newblock Inverse problems and {{Carleman}} estimates.
\newblock {\em Inverse Problems}, 8(4):575--596, 1992.

\bibitem{liu2020reconstruction}
J.-C. Liu and X.-C. Li.
\newblock Reconstruction algorithms of an inverse source problem for the
  {{Helmholtz}} equation.
\newblock {\em Numerical Algorithms}, 84(3):909--933, 2020.

\bibitem{ren2019imaging}
K.~Ren and Y.~Zhong.
\newblock Imaging point sources in heterogeneous environments.
\newblock {\em Inverse Problems}, 35(12):125003, 2019.

\bibitem{wang2019fourier}
X.~Wang, M.~Song, Y.~Guo, H.~Li, and H.~Liu.
\newblock Fourier method for identifying electromagnetic sources with
  multi-frequency far-field data.
\newblock {\em Journal of Computational and Applied Mathematics}, 358:279--292,
  2019.

\end{thebibliography}
\bibliographystyle{abbrv}
\end{document}